\documentclass[reqno,a4paper,twoside]{amsart}
\usepackage{enumitem}

\setenumerate{label=\textnormal{(\arabic*)}}

\newcommand{\dotfillboldidx}[1]{\dotfill{\bf #1}}
\newcommand{\rmidx}[1]{{\rm #1}}

\allowdisplaybreaks[4]

\usepackage{amsmath,amssymb,dsfont,verbatim,bm,geometry,mathrsfs}

\usepackage{mathtools}
\usepackage[raggedright]{titlesec}

\titleformat{\chapter}[display]
{\normalfont\huge\bfseries}{\chaptertitlename\\thechapter}{20pt}{\Huge}
\titleformat{\section}
{\normalfont\Large\bfseries\center}{\thesection}{1em}{}
\titleformat{\subsection}
{\normalfont\large\bfseries}{\thesubsection}{1em}{}
\titleformat{\subsubsection}
{\normalfont\normalsize\bfseries}{\thesubsubsection}{1em}{}
\titleformat{\paragraph}[runin]
{\normalfont\normalsize\bfseries}{\theparagraph}{1em}{}
\titleformat{\subparagraph}[runin]
{\normalfont\normalsize\bfseries}{\thesubparagraph}{1em}{}

\titlespacing*{\chapter} {0pt}{50pt}{40pt}
\titlespacing*{\section} {0pt}{3.5ex plus 1ex minus .2ex}{2.3ex plus .2ex}
\titlespacing*{\subsection} {0pt}{3.25ex plus 1ex minus .2ex}{1.5ex plus .2ex}
\titlespacing*{\subsubsection}{0pt}{3.25ex plus 1ex minus .2ex}{1.5ex plus .2ex}
\titlespacing*{\paragraph} {0pt}{3.25ex plus 1ex minus .2ex}{1em}
\titlespacing*{\subparagraph} {\parindent}{3.25ex plus 1ex minus .2ex}{1em}

\usepackage{tikz}
\usetikzlibrary{matrix, shapes}

\usepackage{float, subfig}
\usepackage{amsrefs}

\renewcommand{\indexname}{Index of symbols}

\usepackage{robustindex}
\usepackage{titletoc}

\usepackage{stackengine}
\usepackage{scalerel}

\contentsmargin{2.55em}
\dottedcontents{section}[1.5em]{}{2em}{1pc}
\dottedcontents{subsection}[4.35em]{}{2.8em}{1pc}
\dottedcontents{subsubsection}[7.6em]{}{3.2em}{1pc}
\dottedcontents{paragraph}[10.3em]{}{3.2em}{1pc}

\usepackage[columns=3,rule=0pt]{idxlayout}
\setindexprenote{For convenience of the reader we list below almost all symbols used in this work, indicating the page in which each of them is introduced. We only except those that are extremely standard.}

\makeindex

\newtheorem{theorem}{Theorem}[section]
\newtheorem{lemma}[theorem]{Lemma}
\newtheorem{proposition}[theorem]{Proposition}
\newtheorem{corollary}[theorem]{Corollary}

\theoremstyle{definition}
\newtheorem{definition}[theorem]{Definition}

\newtheorem{notation}[theorem]{Notation}

\theoremstyle{remark}
\newtheorem{remark}[theorem]{Remark}

\DeclareMathOperator{\Aut}{Aut}
\DeclareMathOperator{\cop}{cop}
\DeclareMathOperator{\codom}{codom}
\DeclareMathOperator{\dom}{dom}
\DeclareMathOperator{\CH}{CH}
\DeclareMathOperator{\coH}{coH}
\DeclareMathOperator{\ide}{id}
\DeclareMathOperator{\Ext}{Ext}
\DeclareMathOperator{\GL}{GL}
\DeclareMathOperator{\Supp}{Supp}
\DeclareMathOperator{\Z}{Z}
\DeclareMathOperator{\en}{en}
\DeclareMathOperator{\En}{En}
\DeclareMathOperator{\factors}{factors}
\DeclareMathOperator{\Ob}{Ob}
\DeclareMathOperator{\op}{op}
\DeclareMathOperator{\Ss}{S}
\DeclareMathOperator{\st}{st}
\DeclareMathOperator{\St}{St}
\DeclareMathOperator{\Cone}{Cone}
\DeclareMathOperator{\PE}{PE}
\DeclareMathOperator{\val}{val}
\DeclareMathOperator{\Valsup}{Valsup}
\DeclareMathOperator{\Valinf}{Valinf}
\DeclareMathOperator{\Succ}{Succ}
\DeclareMathOperator{\Pred}{Pred}
\DeclareMathOperator{\Val}{Val}
\DeclareMathOperator{\ima}{Im}
\DeclareMathOperator{\BC}{BC}
\DeclareMathOperator{\BN}{BN}
\DeclareMathOperator{\BP}{BP}
\DeclareMathOperator{\Ho}{H}
\DeclareMathOperator{\HH}{HH}
\DeclareMathOperator{\HC}{HC}
\DeclareMathOperator{\HP}{HP}
\DeclareMathOperator{\HN}{HN}
\DeclareMathOperator{\Hom}{Hom}
\DeclareMathOperator{\Sh}{Sh}
\DeclareMathOperator{\Tot}{Tot}
\DeclareMathOperator{\Tor}{Tor}
\DeclareMathOperator{\Reg}{Reg}

\newcommand{\fl}{\mathfrak l}
\newcommand{\fr}{\mathfrak r}

\newcommand{\ov}{\overline}
\newcommand{\ot}{\otimes}
\newcommand{\wh}{\widehat}
\newcommand{\wt}{\widetilde}
\newcommand{\ep}{\epsilon}
\newcommand{\De}{\Delta}
\newcommand{\ba}{\mathbf a}
\newcommand{\bx}{\mathbf x}
\newcommand{\bc}{\mathbf c}
\newcommand{\be}{\mathbf e}
\newcommand{\bh}{\mathbf h}
\newcommand{\bv}{\mathbf v}
\newcommand{\byy}{\mathbf y}
\newcommand{\bz}{\mathbf z}
\newcommand{\hs}{\hspace{-0.7pt}}
\newcommand{\hsm}{\hspace{-0.5pt}}
\newcommand{\xcdot}{\hsm\cdot\hsm}
\newcommand{\xcirc}{\hsm\circ\hsm}
\newcommand{\xst}{\hsm*\hsm}

\newcommand{\hatb}{\bar}

\newcommand{\ovb}{\bar}

\newcommand{\wtb}{\bar}

\newcommand{\hJ}{\hs J}

\renewcommand{\theequation}{\thesection.\arabic{equation}}

\DeclareMathAlphabet{\mathpzc}{OT1}{pzc}{m}{it}

\begin{document}

\title[Cleft extensions of weak Hopf algebras]{Cleft extensions of weak Hopf algebras}

\author{Jorge A. Guccione}
\address{Departamento de Matem\'atica\\ Facultad de Ciencias Exactas y Naturales-UBA, Pabell\'on~1-Ciudad Universitaria\\ Intendente Guiraldes 2160 (C1428EGA) Buenos Aires, Argentina.}
\address{Instituto de Investigaciones Matem\'aticas ``Luis A. Santal\'o"\\ Facultad de Ciencias Exactas y Natu\-ra\-les-UBA, Pabell\'on~1-Ciudad Universitaria\\ Intendente Guiraldes 2160 (C1428EGA) Buenos Aires, Argentina.}
\email{vander@dm.uba.ar}

\author{Juan J. Guccione}
\address{Departamento de Matem\'atica\\ Facultad de Ciencias Exactas y Naturales-UBA\\ Pabell\'on~1-Ciudad Universitaria\\ Intendente Guiraldes 2160 (C1428EGA) Buenos Aires, Argentina.}
\address{Instituto Argentino de Matem\'atica-CONICET\\ Savedra 15 3er piso\\ (C1083ACA) Buenos Aires, Argentina.}
\email{jjgucci@dm.uba.ar}

\author[C. Valqui]{Christian Valqui}
\address{Pontificia Universidad Cat\'olica del Per\'u - Instituto de Matem\'atica y Ciencias Afi\-nes, Secci\'on Matem\'aticas, PUCP, Av. Universitaria 1801, San Miguel, Lima 32, Per\'u.}
\email{cvalqui@pucp.edu.pe}


\begin{abstract}
In this paper we study the theory of cleft extensions for a weak bialgebra~$H$. Among other results, we determine when two unitary crossed products of an algebra $A$ by $H$ are equivalent and we prove that if $H$ is a weak Hopf algebra, then the categories of $H$-cleft extensions of an algebra $A$, and of unitary crossed products of $A$ by $H$, are equivalent.
\end{abstract}

\maketitle

\tableofcontents

\section*{Introduction}
The initial motivation for this work comes from our study of the (co)homology of crossed products in weak contexts started in \cite{GGV}, and the present paper provides a theoretical basis for that study. A usual technique for calculating the Hochschild and cyclic homologies of crossed products is to build mixed complexes simpler than the canonical one, which give these homologies. It is convenient that these complexes are provided with filtrations whose associated spectral sequences generalize the classic ones of Hochschild-Serre and Feigin and Tygan. Such complexes exist under very mild conditions. For instance in \cite{CGGV} appropriated complexes were obtained for the Brzezi\'nski's crossed products introduced in \cite{Br}, and it is possible to show that these results remain valid in the more general context considered in~\cites{AFGR, AFGR2, FGR2}. However, if one wants to compute the second face of these spectral sequences it is necessary to restrict ourselves to the crossed products of algebras by weak Hopf algebras with invertible cocycle (or, which is equivalent, to the cleft extensions). For crossed products of algebras with Hopf algebras there is a classical notion of invertible cocycle which works very well. For the crossed products of algebras by weak Hopf algebras introduced in \cite{AFGR} this notion is not so clear, but while we studied them we became convinced that for these crossed products the correct definition of invertible cocycle and of cleft extension are the ones introduced in \cite{AFG}. In that paper and in \cite{AVR} the authors assume that the weak Hopf algebras are cocommutative. This is not a real restriction if one is interested in generalizations of the Sweedler cohomology, but it is  if one is interested in the Hochschild and cyclic Homology of weak crossed products.

Motivated in part by the previous discussion in this paper we consider a symmetric category~$\mathscr{C}$ with split idempotents and we study the weak crossed products with invertible cocycle of an algebra $A$ by a weak Hopf algebra $H$ in $\mathscr{C}$.

\smallskip

Our terminology differs slightly from that used in \cites{AFG, AFGLV, AFGR, AFGR2, AVR, FGR, Ra}. The main differences are two. First that, opposite to the made out in that papers, in the definition of crossed product system we assume that $\nabla_{\!\chi}\xcirc \mathcal{F}=\mathcal{F}$  (this allows us simplify some arguments and it is free of cost because can be achieved simply by replacing $\mathcal{F}$ with $\nabla_{\!\chi}\xcirc \mathcal{F}$); and second, that for us, the crossed product associated with a crossed product system with preunit $(A,V,\chi,\mathcal{F},\nu)$, in which $\mathcal{F}$ is a cocycle satisfying the twisted module condition, is $A\times_{\chi}^{\mathcal{F}} V$ instead of $A\ot _{\chi}^{\mathcal{F}} V$ (but both points of view are equivalent because of Theorem~\ref{prodcruz1}(9)).

\smallskip

This paper is organized as follows: In Section~\ref{Sec: Preliminaries} first we recall the basic properties of the preunits of an algebra; then we review the basic properties of the very general notion of crossed products of algebras by objects of $\mathscr{C}$ introduced in~\cite{AFGR}; and then we do a quick review of the notion of Weak Hopf Algebra. In Section~\ref{subsection: Crossed products by weak Hopf algebras} we begin the study of the crossed products of algebras by weak Hopf algebras introduced in~\cite{AFGR}. The main results are Propositions~\ref{fundamental'},~\ref{caso PiL}, \ref{auxiliar 3} and~\ref{multiplicacion y gamma}, and Corollary~\ref{gamma xst (gamma xcirc Pi^R) = gamma}. In Section~\ref{Sec: Equivalence of weak crossed products} we extend to the setting of arbitrary weak bialgebras the concept of equiv\-a\-lence of crossed products introduced in \cite{AVR}*{Section~3} for the case of cocommutative weak Hopf algebras. Opposite to the made out in that paper we do not require that the cocyle be invertible. In Section~\ref{subsection: Weak crossed products of weak module algebras by weak bialgebras in which the unit cocommutes} we continue the study started in Section~\ref{subsection: Crossed products by weak Hopf algebras}, investigating the consequences that $A$ is a weak $H$-module algebra. The main results are Propositions~\ref{el egalizador es A},~\ref{auxiliar4} and~\ref{conmutatividad de gamma e i}. In Section~\ref{cociclo inversible}, under the assumption that $H$ is a weak Hopf algebra, we prove that each crossed product of a weak $H$-module algebra $A$ by $H$ with invertible cocycle, is an $H$-cleft extension of $A$. Finally, in Section~\ref{Section: Cleft extensions are weak crossed products with invertible cocycle} we prove that each $H$-cleft extension is isomorphic to a crossed product with invertible cocycle, of $H$ by a weak $H$-module algebra. Combining the results of this and the previous section, we obtain that, for each weak Hopf algebra $H$ and each algebra $A$, the categories of unitary crossed products of $A$ by $H$ with invertible cocycle, such that $A$ is a weak $H$-module algebra, is equivalent to the category of $H$-cleft extensions of $A$.

\smallskip

The authors wish to express their gratitude to Jos\'e Nicanor Alonso \'Alvarez y Ram\'on Gonz\'alez Rodr\'iguez for a careful reading of a first version of this paper and several helpful suggestions.

\section{Preliminaries}\label{Sec: Preliminaries}
Throughout this paper  $\mathscr{C}=(\mathcal{C},\ot,K,\alpha,\lambda_l,\lambda_r,c)$ is a symmetric category with split idempotents. We use the words arrow, map and morphism as synonyms. Given objects $U$, $V$ and  $W$ in $\mathscr{C}$ and a map $g\colon V\to W$, we write $U\ot g$ for $\ide_U\ot g$ and $g\ot U$ for $g\ot \ide_U$. By our assumptions, for each idempotent morphism $\phi\colon V \to V$ in $\mathcal{C}$ there exist and object $\phi(V)$ and maps $\imath_{\phi}\colon \phi(V)\to V$ and $p_{\phi}\colon V\to \phi(V)$, that we fix once and for all, such that
$p_{\phi}\circ \imath_{\phi}=\phi$ and $\imath_{\phi}\circ p_{\phi}=\ide_V$.
By the coherence Mac Lane Theorem it suffices to prove the results when $\mathscr{C}$ is strict. In other words, we can assume without loose of generality that the associativity constraint $\alpha$ and the unit constraints $\lambda_l$ and $\lambda_r$ are identities, and we do it. We assume that the reader is familiar with the notions of algebra, coalgebra, module and comodule in $\mathscr{C}$. We also assume that the algebras are associative and the coalgebras are coassociative. Given a unitary algebra $A$ and a counitary coalgebra $C$, we let $\mu\colon A\ot A\to A$, $\eta\colon K\to A$, $\Delta\colon C\to C\ot C$ and $\ep\colon C\to K$
denote the multiplication, the unit, the comultiplication and the counit, respectively, specified with a subscript if necessary.

\smallskip

Given an algebra $A$ and a coalgebra $C$ we consider $\Hom(C,A)$ endowed with the convolution product $\alpha * \beta \coloneqq \mu \xcirc (\alpha\ot \beta)\xcirc \Delta$. It is well known that $\Hom(C,A)$ is a monoid with unit $\eta\xcirc \epsilon$.

\smallskip

We use the nowadays well known graphic calculus for symmetric categories. As usual, morphisms will be composed from top to bottom and tensor products will be represented by horizontal concatenation from left to right. The identity map of an object will be represented by a vertical line. Given an algebra $A$ and a coalgebra $C$, the diagrams
\begin{equation*}\label{eq1}

\end{equation*}
stand for the multiplication map of $A$, the unit of $A$, the action of $A$ on a left $A$-module,the action of $A$ on a right $A$-module, the comultiplication map of $C$, the counit of $C$, the coaction of $C$ on a left $C$-comodule and the coaction of $C$ on a right $C$-comodule, respectively. The natural isomorphism $c$ will be represented by the diagram
\begin{equation*}
%
\end{equation*}
and a map from a tensor product $V_1\ot\cdots \ot V_m$ of $m$ objects to a tensor product $W_1\ot  \cdots \ot W_n$ of~$n$ objects, will be represented by a rectangular box with rounded edges, having $m$ entries up and $n$ outputs down and the name of the arrow inside. We will make the following exceptions: for a twisted space $(A,V,\chi)$ and a map $f\colon H\ot H\to A$, we will use the diagrams
\begin{equation*}
%
\end{equation*}
to represent the map $\chi$, the map $f$ and its convolution inverse, when $f$ is convolution invertible.

\subsection{Preunits}

In this subsection we give a quick review (without proofs) of the notion of preunit introduced in~\cite{CG} (see also \cite{FGR}) and its basic properties.

\begin{definition} Let $\nabla\colon \mathcal{B}\to \mathcal{B}$ be an idempotent morphism. We say that an associative product $\mu_{\mathcal{B}}\colon \mathcal{B}\ot \mathcal{B}\rightarrow \mathcal{B}$ is {\em normalized with respect to $\nabla$} if $\nabla\xcirc\mu_{\mathcal{B}}=\mu_{\mathcal{B}}=\mu_{\mathcal{B}}\xcirc(\nabla\ot \nabla)$.
\end{definition}

\begin{definition}\label{preunit} Let $\mathcal{B}$ be an associative algebra. A {\em preunit of $\mu_{\mathcal{B}}$} is a morphism $\nu\colon K\to \mathcal{B}$ such that $\mu_{\mathcal{B}}\xcirc (\mathcal{B} \ot \nu) = \mu_{\mathcal{B}}\xcirc (\nu \ot \mathcal{B})$ and $\nu=\mu\xcirc (\nu\ot \nu)$.
\end{definition}

Note that if $k$ is a field and $\mathcal{C}$ is the category of $k$-vector spaces, then $\nu\colon K\to \mathcal{B}$ is a preunit of $\mu_{\mathcal{B}}$ if and only if $\nu(1)$ is a central idempotent of $\mathcal{B}$.

\begin{remark}\label{idempotente asociado a una preunidad} Let $\mathcal{B}$ be an associative algebra. If $\nu$ is a preunit of $\mu_{\mathcal{B}}$, then the morphism $\nabla_{\!\nu}\colon \mathcal{B}\to \mathcal{B}$, defined by $\nabla_{\!\nu}\coloneqq \mu_{\mathcal{B}}\xcirc (\mathcal{B}\ot \nu)$, is an idempotent that satisfies
$$
\mu_{\mathcal{B}}\xcirc(\nabla_{\!\nu}\ot \mathcal{B})=\nabla_{\!\nu}\xcirc\mu_{\mathcal{B}}= \mu_{\mathcal{B}}\xcirc(\mathcal{B}\ot \nabla_{\!\nu}) = \mu_{\mathcal{B}}\xcirc(\nabla_{\!\nu}\ot \nabla_{\!\nu}) \qquad\text{and} \qquad \nabla_{\!\nu}\xcirc \nu = \nu.
$$
From this it follows that the associative product $\mu_{\mathcal{B}}^{\nu}\colon \mathcal{B}\ot \mathcal{B}\rightarrow \mathcal{B}$, defined by $\mu_{\mathcal{B}}^{\nu}\coloneqq \nabla_{\!\nu}\xcirc\mu_{\mathcal{B}}$, is normalized with respect to $\nabla_{\!\nu}$. Moreover $\nu$ is a preunit of $\mu_{\mathcal{B}}^{\nu}$ and $\mu_{\mathcal{B}}^{\nu}\xcirc(\mathcal{B}\ot\nu)= \nabla_{\!\nu}$.
\end{remark}

\begin{remark}\label{algebra con unidad asociada a una con preunidad}
Let $\mathcal{B}$ be an algebra and let $\nu$ be a preunit of $\mu_{\mathcal{B}}$. Write $B\coloneqq \nabla_{\!\nu}(\mathcal{B})$, $\imath_{\nu}\coloneqq \imath_{{}_{\nabla_{\!\nu}}}$ and $p_{\nu}\coloneqq p_{{}_{\nabla_{\!\nu}}}$. The map $\mu_{B}\colon B\ot B\rightarrow B$, given by $\mu_{B}\coloneqq p_{\nu}\xcirc\mu_{\mathcal{B}}\xcirc(\imath_{\nu}\ot \imath_{\nu})$, is an associative product with unit $\eta_{B}\coloneqq p_{\nu}\xcirc \nu$. Moreover, $\mu_{B} = p_{\nu}\xcirc\mu_{\mathcal{B}}^{\nu}\xcirc(\imath_{\nu}\ot \imath_{\nu})$, where $\mu_{\mathcal{B}}^{\nu}$ is as in Re\-mark~\ref{idempotente asociado a una preunidad}; the map $p_{\nu}\colon \mathcal{B}\to B$ is multiplicative and the map $\imath_{\nu}\colon B\to \mathcal{B}$ satisfies the equality $\imath_{\nu}\xcirc \mu_{B} = \mu_{\mathcal{B}}^{\nu} \xcirc (\imath_{\nu}\ot \imath_{\nu})$.
\end{remark}

\begin{remark}\label{producto asociado con preunidad}
Let $\nabla\colon \mathcal{B}\to \mathcal{B}$ be an idempotent arrow and let $B\coloneqq \nabla(\mathcal{B})$. If $\mu_{B}\colon B\ot B\rightarrow B$ is an associative product with unit $\eta_{B}$, then the associative product 
$\mu_{\mathcal{B}}\coloneqq \imath_{\nabla}\xcirc\mu_{B}\xcirc(p_{\nabla}\ot p_{\nabla})$ is normalized with respect to $\nabla$ and $\nu\coloneqq \imath_{\nabla}\xcirc \eta_{B}$ is a preunit of $\mu_{\mathcal{B}}$ such that $\nabla  = \nabla_{\!\nu}$.
\end{remark}

\subsection{General weak crossed products}\label{General weak crossed products}

In this subsection we recall a very general notion of crossed product developed in~\cite{AFGR} and~\cite{FGR}, and we review its basic properties.

\begin{definition}\label{no unital twisted space} A triple $(A,V,\chi)$, consisting of an associative unitary algebra $A$ in $\mathcal{C}$, an object~$V$ of $\mathcal{C}$ and a morphism $\chi\colon V\ot A\longrightarrow A\ot V$, is a {\em twisted space} if
\begin{equation}\label{eqtwistingcond}
\chi\xcirc (V\ot \mu_A)= (\mu_A\ot V)\xcirc (A\ot \chi)\xcirc (\chi\ot A).
\end{equation}
In such a case we say that $\chi$ is a {\em twisting map}.
\end{definition}

From here to Definition~\ref{twisted module and cociclo condiciones} inclusive, we assume that $(A,V,\chi)$ is a twisted space. Note that $A\ot V$ is a non unitary $A$-bimodule in $\mathcal{C}$ via the left and right actions $\rho_l$ and $\rho_r$ given by
\begin{equation}\label{etq}
\rho_l\coloneqq\mu_A\ot V\qquad\text{and}\qquad \rho_r\coloneqq(\mu_A\ot V)\xcirc (A\ot \chi),
\end{equation}
respectively. Clearly $\rho_l$ is an unitary action. We let $\nabla_{\!\chi}$ denote the left and right $A$-linear idem\-potent en\-do\-morph\-ism of $A\ot V$, defined by $\nabla_{\!\chi}\coloneqq \rho_r\xcirc (A\ot V\ot \eta_A)$. Set $A\times V\coloneqq \nabla_{\!\chi}(A\ot V)$, $\imath_{\chi}\coloneqq \imath_{_{\nabla_{\!\chi}}}$ and $p_{\chi}\coloneqq p_{_{\nabla_{\!\chi}}}$. Note that $A\times V$ is an unitary $A$-bimodule via the left and right actions $\bar{\rho}_l$ and $\bar{\rho}_r$ given by $\bar{\rho}_l \!\coloneqq\! p_{\chi}\xcirc \rho_l\xcirc (A\ot \imath_{\chi})$ and $\bar{\rho}_r\!\coloneqq\! p_{\chi}\xcirc \rho_r\xcirc (\imath_{\chi}\ot A)$, re\-specti\-vely. Moreover, $\imath_{\chi}$ and $p_{\chi}$ are $A$-bimodule morphisms.

\begin{definition}\label{crossed product system} A tuple $(A,V,\chi,\mathcal{F})$ is a {\em crossed product system} if $\mathcal{F}\colon V\ot V\longrightarrow A\ot V$ is a morphism such that $\nabla_{\!\chi}\xcirc \mathcal{F}=\mathcal{F}$.
\end{definition}

\begin{definition}[\cite{AFGR}*{Definitions 2.4 and 2.5}]\label{twisted module and cociclo condiciones} We say that $\mathcal{F}$ is a {\em cocycle} that satisfies the {\em twisted module condition} if
\begin{equation*}

\end{equation*}
\noindent More precisely, the first equality says that $\mathcal{F}$  satisfies the twisted module condition and the second one says that $\mathcal{F}$ is a cocycle.
\end{definition}

From here to Definition~\ref{crossed product system with preunit} inclusive $(A,V,\chi,\mathcal{F})$ is a crossed product system.

\begin{notation}\label{mult A otimes V} We let $A\ot _{\chi}^{\mathcal{F}} V$ and $A\times_{\chi}^{\mathcal{F}} V$ denote the objects $A\ot V$ and $A\times V$, endowed with the multiplication maps $\mu_{A\ot _{\chi}^{\mathcal{F}} V}$ and $\mu_{A\times_{\chi}^{\mathcal{F}} V}$, defined by
$$
\mu_{A\ot _{\chi}^{\mathcal{F}} V}\coloneqq (\mu_A\ot V)\xcirc(\mu_A\ot \mathcal{F}) \xcirc(A\ot \chi\ot V)\quad\text{and}\quad
\mu_{A\times_{\chi}^{\mathcal{F}} V}\coloneqq p_{\chi}\xcirc\mu_{A\ot _{\chi}^{\mathcal{F}} V} \xcirc(\imath_{\chi}\ot \imath_{\chi}),
$$
respectively. For the sake of brevity we will write $E\coloneqq$ instead of $\cramped{A\times_{\chi}^{\mathcal{F}}} V$ and $\mathcal{E}$ instead of $\cramped{A\ot _{\chi}^{\mathcal{F}}} V$. We consider $\mathcal{E}\ot \mathcal{E}$ and $E\ot E$ as $A$-bimodules in the natural way.
\end{notation}

\begin{definition}\label{producto cruzado} We say that $E$ is the {\em crossed product of $A$ by $V$ associated with $\chi$ and $\mathcal{F}$} if~$\mathcal{F}$ a cocycle that satisfies the twisted module condition.
\end{definition}

\begin{definition}\label{crossed product system with preunit} Let $\nu\colon K\rightarrow A\ot V$ be an arrow. The tuple $(A,V,\chi,\mathcal{F},\nu)$ is a {\em crossed product system with preunit} if
\begin{align}
& (\mu_A\ot V)\xcirc(A\ot \mathcal{F})\xcirc(\chi\ot V)\xcirc(V\ot \nu) = \nabla_{\!\chi}\xcirc (\eta_A\ot V),\label{preunit1}\\
& (\mu_A\ot V)\xcirc (A\ot \mathcal{F})\xcirc(\nu\ot V) = \nabla_{\!\chi}\xcirc (\eta_A\ot V),\label{preunit2}\\
& (\mu_A\ot V)\xcirc (A\ot \chi)\xcirc(\nu\ot A) = (\mu_A\ot V)\xcirc (A\ot \nu).\label{preunit3}
\end{align}
\end{definition}

Let $(A,V,\chi,\mathcal{F},\nu)$ be a crossed product system with preunit and let
$$
\nabla_{\!\nu}\colon A\ot V\longrightarrow A\ot V,\quad \jmath'_{\nu}\colon A \to \mathcal{E},\quad \jmath_{\nu}\colon A\to E \quad\text{and}\quad \gamma\colon V\to E
$$
be the arrows defined by
$$
\nabla_{\!\nu}\coloneqq \mu_\mathcal{{E}}\xcirc \bigl(\mathcal{E} \ot \nu\bigr),\quad \jmath'_{\nu}\coloneqq (\mu_A\ot V)\xcirc (A\ot \nu),\quad \jmath_{\nu}\coloneqq p_{\chi}\xcirc \jmath'_{\nu}\quad\text{and}\quad \gamma\coloneqq p_{\chi}\xcirc (\eta_A \ot V).
$$

\smallskip

Except for items 4), ~7) and the assertions about the right $A$-linearity in items 3), 5) and~6), whose proofs we leave to the reader, the following result is \cite{FGR}*{Remark~3.10, one implication of Theorem~3.11 and Corollary~3.12}.

\begin{theorem}\label{prodcruz1} Let $(A,V,\chi,\mathcal{F},\nu)$ be a crossed product system with preunit. If $\mathcal{F}$ is a cocycle that satisfies the twisted module condition, then the following facts hold:

\begin{enumerate}[itemsep=0.7ex,  topsep=1.0ex, label=\emph{(\arabic*)}]

\item $\mu_{\mathcal{E}}$ is a left and right $A$-linear associative product, that is normalized with respect to~$\nabla_{\!\chi}$.

\item $\nu$ is a preunit of $\mu_{\mathcal{E}}$, $\nabla_{\!\nu}\xcirc\nu = \nu$ and $\nabla_{\!\nu} = \nabla_{\!\chi}$.

\item $\mu_E$ is left and right $A$-linear, associative and has unit $\eta_E\coloneqq p_{\chi}\xcirc \nu$.

\item The maps $\imath_{\chi}$ and $p_{\chi}$ are multiplicative.

\item $\jmath'_{\nu}$ is left and right $A$-linear, multiplicative, and satisfies $\nabla_{\!\nu} \xcirc \jmath'_{\nu} = \jmath'_{\nu}$.

\item $\jmath_{\nu}$ is left and right $A$-linear, multiplicative and unitary.

\item $\mu_E\xcirc (\jmath_{\nu}\ot E) = \bar{\rho}_l$ and $ \mu_E\xcirc (E \ot \jmath_{\nu}) = \bar{\rho}_r$.

\item $\chi = \mu_{\mathcal{E}}\xcirc (\eta\ot V \ot \jmath'_{\nu})$ and $\mathcal{F} = \mu_{\mathcal{E}}\xcirc (\eta\ot V \ot \eta\ot V)$.

\item $\chi = \imath_{\chi}\xcirc\mu_E \xcirc ( \gamma \ot \jmath_{\nu})$ and $\mathcal{F} = \imath_{\chi}\xcirc \mu_E\xcirc (\gamma \ot \gamma)$.

\end{enumerate}

\end{theorem}

\begin{definition}\label{producto cruzado unitario} Let $(A,V,\chi,\mathcal{F}, \nu)$ be a crossed product system with preunit. If $\mathcal{F}$ a cocycle that satisfies the twisted module condition, then we say that the algebra $E$ is the {\em unitary crossed product of $A$ by $V$ associated with $\chi$, $\mathcal{F}$ and $\nu$}. We also say that $\chi$ and $\mathcal{F}$ are the twisting map and the cocycle of $E$, respectively.
\end{definition}

\begin{remark}\label{propiedad de gamma} By Theorem~\ref{prodcruz1}(7) and the fact that $p_{\!\chi}$ is left $A$-linear and $\gamma =p_{\chi}\xcirc (\eta_A \ot V)$,
\begin{equation}\label{iota(a)gamma(v)}
\mu_E\xcirc (\jmath_{\nu}\ot \gamma) = \bar{\rho}_l\xcirc (A\ot \gamma) = p_{\chi}.
\end{equation}
\end{remark}

\begin{remark}\label{gama iota y gama gama} For each crossed product system with preunit $(A,V,\chi,\mathcal{F},\nu)$, we have
$$
\mu_E\xcirc (\gamma \ot \jmath_{\nu}) = \mu_E\xcirc (\jmath_{\nu}\ot \gamma) \xcirc \chi\quad\text{and}\quad \mu_E\xcirc (\gamma\ot \gamma)= \mu_E\xcirc (\jmath_{\nu}\ot \gamma) \xcirc \mathcal{F}.
$$
\end{remark}

\subsection{Weak Hopf Algebras}\label{subsection: Weak Hopf algebras}
Weak bialgebras and weak Hopf algebras are generalizations of bialgebras and Hopf algebras, introduced in~\cites{BNS1,BNS2}, in which the axioms about the unit, the counit and the antipode are replaced by weaker properties. For the convenience of the reader, in this subsection we collet without proofs the properties of weak Hopf algebras in an symmetric tensor category $\mathcal{C}$, with split idempotents, that we need in this paper. All the results considered by us were established in~\cites{BNS1, CG}, or they are immediate consequence of results obtained in those papers. In spite of that in~\cites{BNS1, CG} the authors work in the setting of finite dimensional vector spaces, all the results in this subsection (and in this paper) are valid in the context of symmetric tensor categories.

\begin{definition}\label{weak bialgebra} A {\em weak bialgebra in $\mathcal{C}$} is an object $H$ endowed with an unitary algebra and a counitary coalgebra structure, such that:

\begin{enumerate}[itemsep=0.7ex,  topsep=1.0ex, label=(\arabic*)]

\item $\Delta \xcirc \mu = (\mu\ot \mu) \xcirc c\xcirc (\Delta \ot \Delta)$,

\item $\epsilon\xcirc \mu\xcirc (\mu\ot H)= (\epsilon \xcirc \mu \ot \epsilon \xcirc \mu)\xcirc (H\ot \Delta \ot H)= (\epsilon \xcirc \mu \ot \epsilon \xcirc \mu)\xcirc (H\ot \Delta^{\cop} \ot H)$,

\item $(\Delta \ot H)\xcirc \Delta \xcirc \eta= (H\ot \mu \ot H)\xcirc (\Delta\xcirc \eta \ot \Delta\xcirc \eta)=(H\ot \mu^{\op} \ot H)\xcirc (\Delta\xcirc \eta \ot \Delta\xcirc \eta)$,

\end{enumerate}
where $\Delta^{\cop}\coloneqq c\xcirc \Delta$ and $\mu^{\op}\coloneqq \mu\xcirc c$.
\end{definition}

For a weak bialgebra $H$, we denote by $\Pi^{\hs L}$, $\Pi^{\hs R}$, $\ov{\Pi}^L$ and $\ov{\Pi}^R$ the maps, defined by
\begin{align*}
&\Pi^{\hs L}\coloneqq (\epsilon\xcirc \mu \ot H) \xcirc (H \ot c) \xcirc (\Delta\xcirc\eta \ot H), && \ov{\Pi}^L\coloneqq (H\ot \epsilon \xcirc\mu) \xcirc (\Delta\xcirc\eta\ot H),\\
&\Pi^{\hs R}\coloneqq (H \ot \epsilon\xcirc \mu) \xcirc (c \ot H) \xcirc (H\ot \Delta\xcirc\eta), && \ov{\Pi}^R\coloneqq (\epsilon \xcirc\mu \ot H) \xcirc (H\ot \Delta\xcirc\eta).
\end{align*}
A direct computation shows that $\Pi^{\hs L}$, $\Pi^{\hs R}$, $\ov{\Pi}^{\hs L}$ and $\ov{\Pi}^{\hs R}$ are unitary and counitary idempotent morphisms (see, for instance, \cites{BNS1,CG}). For each $X\in\{L,R\}$, set $H^{\hs X}\coloneqq \Pi^{\hs X}(H)$, $\imath_{\hs X}\coloneqq \imath_{\Pi^{\!X}}$ and $p_{\hs X}\coloneqq p_{\Pi^{\!X}}$.

\begin{proposition}\label{algunas composiciones} The equalities $\ov{\Pi}^L\xcirc \Pi^R = \Pi^R$, $\ov{\Pi}^R\xcirc \Pi^L = \Pi^L$ and $\Pi^R\xcirc \ov{\Pi}^L = \ov{\Pi}^L$ hold.
%
\end{proposition}


\begin{remark}\label{ide ast Pi^R = ide and} An immediate computation shows that $\ide*\Pi^{\hs R}=\ide = \Pi^{\hs L}*\ide$.
\end{remark}


\begin{proposition}\label{Delta compuesto con unidad incluido en H^{hs R}ot H^{hs L}} The following equalities hold:

\begin{enumerate}[itemsep=0.7ex,  topsep=1.0ex, label=\emph{(\arabic*)}]


\item $\Delta\xcirc \eta = (H\ot \Pi^{\hs L})\xcirc \Delta\xcirc \eta = (\Pi^{\hs R}\ot H)\xcirc \Delta\xcirc \eta = (\Pi^{\hs R}\ot \Pi^{\hs L})\xcirc \Delta\xcirc \eta$,


\item $\Delta\xcirc \eta = (H\ot \ov{\Pi}^R)\xcirc \Delta\xcirc \eta = (\ov{\Pi}^L\ot H)\xcirc \Delta\xcirc \eta = (\ov{\Pi}^L\ot \ov{\Pi}^R)\xcirc \Delta\xcirc \eta$.

\end{enumerate}

\end{proposition}

\begin{proposition}\label{subalgebras} The following equalities hold:

\begin{enumerate}[itemsep=0.7ex,  topsep=1.0ex, label=\emph{(\arabic*)}]

\item $\mu\xcirc (\Pi^{\hs L}\ot\Pi^{\hs L}) = \Pi^{\hs L}\xcirc \mu\xcirc (\Pi^{\hs L}\ot\Pi^{\hs L}) = \Pi^{\hs L}\xcirc \mu\xcirc (\Pi^{\hs L}\ot H)$,

\item $\mu\xcirc (\Pi^{\hs R}\ot\Pi^{\hs R}) = \Pi^{\hs R}\xcirc \mu\xcirc (\Pi^{\hs R}\ot\Pi^{\hs R}) = \Pi^{\hs R}\xcirc \mu\xcirc (H\ot\Pi^{\hs R})$.



\end{enumerate}
\end{proposition}

\begin{remark}\label{HR y HL son subalgebras} From the fact that $\Pi^L$ and $\Pi^R$ are unitary maps and Proposition~\ref{subalgebras} it follows that  $H^{\hs R}$ and $H^{\hs L}$ are unitary associative algebras via the multiplication maps $p_{\hs R}\xcirc \mu \xcirc (\imath_{\hs R}\ot\imath_{\hs R})$ and $p_{\hs L}\xcirc \mu \xcirc (\imath_{\hs L}\ot\imath_{\hs L})$, respectively. Moreover the arrows $\imath_{\hs R}$ and $\imath_{\hs L}$ are algebra morphisms.
\end{remark}

\begin{proposition}\label{conmutatividad debil} It is true that $\mu\xcirc c\xcirc (\Pi^{\hs L}\ot\Pi^{\hs R}) = \mu\xcirc (\Pi^{\hs L}\ot\Pi^{\hs R})$.
\end{proposition}

\begin{proposition}\label{mu Pi^R, etc} The following equalities hold:
\begin{enumerate}[itemsep=0.7ex,  topsep=1.0ex, label=\emph{(\arabic*)}]

\item $\mu\xcirc (H \ot \Pi^{\hs L}) = (\epsilon\xcirc \mu\ot H)\xcirc (H\ot c)\xcirc (\Delta\ot H)$,
\item $\mu\xcirc (\Pi^{\hs R} \ot H) =  (H\ot \epsilon\xcirc \mu)\xcirc (c\ot H)\xcirc (H\ot \Delta)$.
\end{enumerate}
\end{proposition}

\begin{proposition}\label{delta Pi^R, etc} The equality $(H \ot \Pi^{\hs L}) \xcirc \Delta = (\mu\ot H)\xcirc (H\ot c)\xcirc (\Delta\xcirc \eta\ot H)$ holds.
%
\end{proposition}

\begin{proposition}\label{mu delta Pi^R, etc} The following equalities hold:
\begin{align*}
& \Pi^{\hs L} \xcirc \mu = \Pi^{\hs L} \xcirc \mu \xcirc (H\ot \Pi^{\hs L}), && \Delta\xcirc \Pi^{\hs L}  = (H\ot \Pi^{\hs L})\xcirc \Delta\xcirc \Pi^{\hs L},\\
& \Pi^{\hs R} \xcirc \mu = \Pi^{\hs R} \xcirc \mu \xcirc (\Pi^{\hs R}\ot H), && \Delta\xcirc \Pi^{\hs R}  = (\Pi^{\hs R} \ot H)\xcirc \Delta\xcirc \Pi^{\hs R}.
\end{align*}
\end{proposition}

\begin{proposition}\label{Delta compuesto con Pi^R} The following equalities hold:

\begin{enumerate}[itemsep=0.7ex,  topsep=1.0ex, label=\emph{(\arabic*)}]

\item $(\Pi^{\hs L}\ot \epsilon\xcirc  \mu)\xcirc(\Delta\ot H)= \Pi^{\hs L}\xcirc \mu =  (\epsilon\xcirc \mu\ot \Pi^{\hs L})\xcirc (H\ot c) \xcirc(\Delta\ot H)$,

\item $ (\epsilon\xcirc\mu\ot \Pi^{\hs R}) \xcirc (H\ot \Delta) = \Pi^{\hs R}\xcirc \mu =  (\Pi^{\hs R}\ot \epsilon\xcirc\mu) \xcirc (c\ot H)\xcirc (H\ot \Delta)$,

\item $(\mu\ot H)\xcirc (\Pi^{\hs L}\ot \Delta\xcirc \eta)=\Delta\xcirc \Pi^{\hs L} = (\mu \ot H) \xcirc (H\ot c)\xcirc (\Delta\xcirc \eta\ot \Pi^{\hs L})$.

\item $(H\ot \mu)\xcirc (\Delta \xcirc \eta\ot \Pi^{\hs R})=\Delta\xcirc \Pi^{\hs R} = (H \ot \mu) \xcirc (c\ot H)\xcirc (\Pi^{\hs R}\ot \Delta\xcirc \eta)$,

\end{enumerate}
\end{proposition}

\begin{proposition}\label{Delta compuesto con Pi^R multiplicado por algo} Let $H$ be a weak bialgebra. Then:
\begin{align*}
& \Delta \xcirc \mu\xcirc (\Pi^{\hs R} \ot H) =   (H\ot \mu)\xcirc (c\ot H) \xcirc (\Pi^{\hs R}\!\ot \Delta), &&\Delta \xcirc \mu\xcirc (H\ot \Pi^{\hs R}) = (H\ot \mu)\xcirc  (\Delta\ot \Pi^{\hs R}),\\
& \Delta \xcirc \mu\xcirc (H\ot \Pi^{\hs L})= (\mu\ot H)\xcirc (H\ot c)\xcirc (\Delta\ot \Pi^{\hs L}), &&\Delta \xcirc \mu\xcirc (\Pi^{\hs L}\ot H)= (\mu\ot H)\xcirc (\Pi^{\hs L}\ot \Delta).
\end{align*}
\end{proposition}

\begin{definition}\label{weak Hopf Algebra} Let $H$ be a weak bialgebra. An {\em antipode} of $H$ is a map $S\colon H\to H$ (or $S_H$ if necessary to avoid confusion), that have the following properties:

\begin{enumerate}[itemsep=0.7ex,  topsep=1.0ex, label=(\arabic*)]

\item $\mu\xcirc (H\ot S)\xcirc \Delta = \Pi^{\hs L}$,

\item $\mu\xcirc (S\ot H)\xcirc \Delta = \Pi^{\hs R}$,

\item $\mu\xcirc (\mu\ot H)\xcirc (S\ot H \ot S)\xcirc (\Delta\ot H)\xcirc \Delta=S$.

\smallskip

\end{enumerate}
We say that $H$ is a {\em weak Hopf algebra} if it has an antipode.
\end{definition}

As it was shown in~\cite{BNS1}, an antipode $S$, if there exists, is unique. It was also shown in~\cite{BNS1} that~$S$ is antimultiplicative, anticomultiplicative and leaves the unit and counit invariant.

\begin{proposition}\label{S y Pi} It is true that $\Pi^{\hs L} = \ov{\Pi}^R\xcirc S = S\xcirc \ov{\Pi}^L$ and $\Pi^{\hs R} = \ov{\Pi}^L\xcirc S = S\xcirc \ov{\Pi}^R$.
\end{proposition}

\begin{proposition}\label{Delta eta con S} Let $H$ be a weak Hopf algebra. The following equalities hold:

\begin{enumerate}[itemsep=1.0ex, topsep=1.0ex, label=\emph{(\arabic*)}]

\item $(\mu\ot S)\xcirc (H\ot \Delta)\xcirc (H\ot \eta) = (H\ot \Pi^R)\xcirc \Delta$,

\item $(S\ot \mu)\xcirc (\Delta \ot H)\xcirc (\eta\ot H) = (\Pi^L\ot H)\xcirc \Delta$.

\end{enumerate}
\end{proposition}

\section{Crossed products by weak bialgebras}\label{subsection: Crossed products by weak Hopf algebras}
In this section we begin the study of the crossed product of an algebra $A$ by a weak bialgebra~$H$, introduced in \cite{AFGR}.

\subsection{Basic properties}

\begin{definition} An arrow $\rho\colon H\ot A\longrightarrow A$ is a {\em weak measure of $H$ on $A$} if
$$
\rho\xcirc (H\ot \mu_A) = \mu \xcirc (\rho \ot \rho) \xcirc (H\ot c\ot A) \xcirc (\Delta \ot A\ot A).
$$
\end{definition}

From here to the end of this section $\rho$ is a weak measure. Let $\chi_{\rho}\colon H\ot A\longrightarrow A\ot H$ be the morphism defined by $\chi_{\rho}\coloneqq (\rho\ot H) \xcirc (H\ot c) \xcirc (\Delta\ot A)$. It is well known that $(A,H,\chi_{\rho})$ is a twisted space. Moreover $(A\ot \epsilon)\xcirc \chi_{\rho} = \rho$. Recall that $A\ot H$ is a non unitary $A$-bimodule via~\eqref{etq}. Let $A\times H$, $\nabla_{\!\chi_{\rho}}$, $p_{\!\chi_{\rho}}$ and $\imath_{\!\chi_{\rho}}$ be as in Subsection~\ref{General weak crossed products}. Set $\nabla_{\!\rho}\coloneqq \nabla_{\!\chi_{\rho}}$, $p_{\rho}\coloneqq p_{\!\chi_{\rho}}$ and $\imath_{\rho}\coloneqq \imath_{\!\chi_{\rho}}$. By definition
\begin{equation}\label{calculo de nabla rho} \nabla_{\!\rho}=  (\mu \ot H) \xcirc (A\ot \rho\ot H) \xcirc (A\ot H\ot c) \xcirc (A \ot \Delta \ot \eta_A).
\end{equation}
Furthermore, we know that $\chi_{\rho}= \nabla_{\!\rho} \xcirc \chi_{\rho}$, that $\nabla_{\!\rho}$ is left and right $A$-linear and idempotent and that $p_{\rho}$ and $\imath_{\rho}$ are left and right $A$-linear.

\smallskip

Given a morphism $f\colon H\ot H\rightarrow A$, we define $\mathcal{F}_{\hs f}\colon H\ot H\longrightarrow A\ot H$ by $\mathcal{F}_{\hs f}\coloneqq (f\ot \mu) \xcirc \Delta_{H^{\ot^2}}$.

\begin{remark}\label{aux} A direct computation using, that $c$ is natural, $\Delta$ and $\Delta_{H^{\ot^2}}$ are coassociative and item~1) of Definition~\ref{weak bialgebra}, shows that
\begin{enumerate}[itemsep=0.7ex,  topsep=1.0ex, label=(\arabic*)]

\item $(A\ot \Delta)\xcirc \chi_{\rho} = (\chi_{\rho}\ot H)\xcirc (H\ot s) \xcirc (\Delta\ot A)$,

\item $(A\ot\Delta)\xcirc \mathcal{F}_{\hs f} = (\mathcal{F}_{\hs f}\ot\mu)\xcirc \Delta_{H^{\ot^2}}$,

\item $(H^{\ot^2}\ot \mathcal{F}_{\hs f})\xcirc \Delta_{H^{\ot^2}} =
(H^{\ot^2}\ot f\ot \mu)\xcirc (\Delta_{H^{\ot^2}}\ot \mu)\xcirc \Delta_{H^{\ot^2}}$.

\end{enumerate}
\end{remark}

\begin{remark}\label{coaccion en A ot H}
It is evident that $A\ot H$ is a counitary $H$-comodule via the map $\delta_{A\ot H}\coloneqq A\ot \Delta$.  Moreover, by Remark~\ref{aux}(1), the map $\nabla_{\!\rho}$ is $H$-colinear. A direct computation using this shows that $A\times H$ is a counitary $H$-comodule via $\delta_{A\times H}\coloneqq (p_{\rho}\ot H) \xcirc \delta_{A\ot H} \xcirc \imath_{\rho}$ and that $\imath_{\rho}$ and $p_{\rho}$ are $H$-colinear maps.
\end{remark}

\begin{proposition}\label{equivalencia de que F incluido en A times V} The following conditions are equivalent:

\begin{enumerate}[itemsep=0.7ex,  topsep=1.0ex, label=\emph{(\arabic*)}]

\item $f = (A\ot \epsilon)\xcirc \mathcal{F}_{\hs f}$ and $\mathcal{F}_{\hs f} = \nabla_{\!\rho} \xcirc \mathcal{F}_{\hs f}$.

\item $f = \mu_A\xcirc(A \ot\rho)\xcirc (\mathcal{F}_{\hs f} \ot \eta_A)$.

\end{enumerate}

\end{proposition}

\begin{proof} 1) $\Rightarrow$ 2)\enspace We have
$
f = (A\ot\epsilon)\xcirc \mathcal{F}_{\hs f} = (A\ot\epsilon)\xcirc \nabla_{\!\rho}\xcirc \mathcal{F}_{\hs f} = \mu_A\xcirc(A \ot\rho)\xcirc (\mathcal{F}_{\hs f} \ot \eta_A).
$

\smallskip

\noindent 2) $\Rightarrow$ 1)\enspace Left to the reader (use Remark~\ref{aux}(2)).
\end{proof}

Consequently, if $f$ satisfies Proposition~\ref{equivalencia de que F incluido en A times V}(2), then $(A,H,\chi_{\rho},\mathcal{F}_{\hs f})$ is a crossed product system. If this is the case, then we set $E= A\times_{\rho}^f V\coloneqq A\times_{\chi_{\rho}}^{\mathcal{F}_{\hs f}} V$ and  $\mathcal{E} = A\ot_{\rho}^f V\coloneqq A\ot_{\chi_{\rho}}^{\mathcal{F}_{\hs f}} V$.

\begin{definition} A map $f\colon H\ot H\to A$ is a cocycle that satisfies the twisted module condition~if
$$

$$
More precisely, the first equality says that $f$ satisfies the twisted module condition and the second one says that $f$ is a cocycle.
\end{definition}

The following result was established in the proof of~\cite{FGR}*{Theorem~4.2}. For the convenience of the reader we provide a diagrammatic proof.

\begin{proposition}\label{cocy equiv cocy y..} If $f = (A\ot \epsilon)\xcirc \mathcal{F}_{\hs f}$, then the following assertions are true:

\begin{enumerate}[itemsep=0.7ex,  topsep=1.0ex, label=\emph{(\arabic*)}]

\item $f$ satisfies the twisted module condition if and only if $\mathcal{F}_{\hs f}$ does it.

\item $f$ is a cocycle if and only if $\mathcal{F}_{\hs f}$ is.

\end{enumerate}
\end{proposition}

\begin{proof} We prove item~(1) and leave the proof of item~(2) to the reader. Composing the first equality of Definition~\ref{twisted module and cociclo condiciones} with $A\ot\epsilon$ and using that $\rho = (A\ot \epsilon)\xcirc \chi_{\rho}$ and $f = (A\ot \epsilon)\xcirc \mathcal{F}_{\hs f}$, we obtain that if $\mathcal{F}_{\hs f}$ satisfies the twisted module condition, then $f$ also satisfies the twisted module condition. Conversely, if this happens, then we have
$$

$$
where the first equality holds by the very definition of $\mathcal{F}_{\hs f}$; the second one, by Remark~\ref{aux}(1); the third one, since $c$ is natural; the fourth one, since $f$ satisfies the twisted module condition; the fifth one, by Remark~\ref{aux}(2); and the last one, by the very definition of $\chi_{\rho}$.
\end{proof}

\begin{proposition}\label{fundamental'} If $f\colon H\ot H\to A$ satisfies $f = (A\ot \epsilon)\xcirc \mathcal{F}_{\hs f}$, then
$$
f \xcirc (\mu\ot H)\xcirc (H\ot \Pi^{\hs R} \ot H) = f\xcirc (H\ot \mu)\xcirc (H\ot \Pi^{\hs R} \ot H).
$$
\end{proposition}

\begin{proof} By the hypothesis, Proposition~\ref{Delta compuesto con Pi^R multiplicado por algo}, the fact that $c$ is natural and $\mu_H$ is associative
$$

$$
as desired.
\end{proof}

\begin{proposition}\label{caso PiL} If $f\colon H\ot H\to A$ satisfies $(\epsilon\xcirc \mu\ot f)\xcirc \Delta_{H^{\ot^2}}=f$, then
$$
f \xcirc (\mu\ot H)\xcirc (H\ot \Pi^L \ot H) = f\xcirc (H\ot \mu)\xcirc (H\ot \Pi^L \ot H).
$$
\end{proposition}

\begin{proof} Mimic the proof of Proposition~\ref{fundamental'}.
\end{proof}

\begin{proposition}[\cite{Ra}*{Remark~3.5}]\label{condicion 1 sobre 1_E} Let $\nu\colon K\to A\ot H$ be a map. The following assertions are equivalent:
\begin{enumerate}[itemsep=0.7ex,  topsep=1.0ex, label=\emph{(\arabic*)}]

\item $\nu = (A\ot\Pi^{\hs L})\xcirc \nu$.

\item $(A\ot\Delta)\xcirc\nu = (A\ot \mu\ot H)\xcirc(A\ot H\ot \Delta) \xcirc(\nu\ot\eta)$.

\item $(A\ot\Delta)\xcirc\nu = (A\ot \mu\ot H)\xcirc(A\ot c\ot H)\xcirc(A\ot H\ot \Delta) \xcirc(\nu\ot\eta)$.

\end{enumerate}

\end{proposition}

\begin{proof} (1) $\Rightarrow$ (2) and~(3)\enspace By Proposition~\ref{Delta compuesto con Pi^R}(3) and the fact that $\Pi^{\hs L}$ is idempotent, we have
$$

$$
as desired.

\smallskip

\noindent (2) $\Rightarrow$ (1)\enspace Since $\Delta$ is counitary this follows composing $A\ot\epsilon \ot H$ to the equality in item~(2) and using Proposition~\ref{Delta compuesto con unidad incluido en H^{hs R}ot H^{hs L}}(1).

\smallskip

\noindent (3) $\Rightarrow$ (1)\enspace Mimic the proof of (2) $\Rightarrow$ (1).
\end{proof}

The following result was established in the proof of~\cite{FGR}*{Theorem~4.2}. For the convenience of the reader we provide a diagrammatic proof.

\begin{proposition}\label{cond1preunudad} Let $\nu\colon K\rightarrow A\ot H$ and $f\colon H\ot H\to A$ be maps. Assume that $f$ satisfies the first equality in Proposition~\ref{equivalencia de que F incluido en A times V}(1) and that $\nu = (A\ot\Pi^L)\xcirc \nu$. The following assertions hold:

\begin{enumerate}[itemsep=0.7ex,  topsep=1.0ex, label=\emph{(\arabic*)}]

\item The map $\nu$ satisfies \eqref{preunit1} if and only if $\rho\xcirc (H\ot \eta_A) = \mu_A\xcirc (\rho\ot f) \xcirc (H\ot c\ot H)\xcirc (\Delta\ot\nu)$.

\item The map $\nu$ satisfies \eqref{preunit2} if and only if $\rho\xcirc (H\ot \eta_A) = \mu_A\xcirc (A\ot f)\xcirc (\nu\ot H)$.

\end{enumerate}

\end{proposition}

\begin{proof} (1)\enspace Assume $\rho\xcirc (H\ot \eta_A) = \mu_A\xcirc (\rho\ot f) \xcirc (H\ot c\ot H)\xcirc (\Delta\ot\nu)$ and let $\daleth\coloneqq (H^{\ot^2}\ot \mu)\xcirc \Delta_{H^{\ot^2}}$. Then, by the coas\-so\-cia\-tiv\-i\-ty of $\Delta$, Definition~\ref{weak bialgebra}(1), the fact that $\mu$ is unitary and $c$ is natural, Propositions~\ref{Delta compuesto con unidad incluido en H^{hs R}ot H^{hs L}}(1) and~\ref{fundamental'} and condition~(3) in Proposition~\ref{condicion 1 sobre 1_E},
$$

$$
which is condition~\eqref{preunit1}. In order to prove the converse, it suffices to apply $A\ot \epsilon$ to the equality in~\eqref{preunit1} and use the first equality Proposition~\ref{equivalencia de que F incluido en A times V}(1).

\smallskip

\noindent (2)\enspace Mimic the proof of item~(1).
\end{proof}

Let $\nu\colon K\rightarrow A\ot H$ be the map $\nu\coloneqq \nabla_{\!\rho}\xcirc (\eta_A\ot \eta)$. Clearly
\begin{equation}\label{cal de gamma}
\nu = \chi_{\!\rho}\xcirc (\eta\ot \eta_A),\quad \nu = (A\ot \Pi^{\hs L})\xcirc \nu \quad\text{and}\quad \imath_{\rho}\xcirc \gamma = \nabla_{\!\rho}\xcirc (\eta_A\ot H) = \chi_{\!\rho}\xcirc (H\ot \eta_A),
\end{equation}
where $\gamma$ is as above of Theorem~\ref{prodcruz1}. Note that
\begin{equation}\label{cal de gamma'}
(A\ot \epsilon)\xcirc \nu = (A\ot \epsilon)\xcirc \chi_{\!\rho}\xcirc (\eta\ot \eta_A) = \rho \xcirc (\eta\ot \eta_A)
\end{equation}
and that by Remark~\ref{coaccion en A ot H} the map $\gamma$ is $H$-colinear.

\begin{theorem}[\cite{Ra}*{Proposition~3.6}]\label{weak crossed prod} Let $f\colon H\ot H\rightarrow A$ be a map.  If

\begin{enumerate}[itemsep=0.7ex,  topsep=1.0ex, label=\emph{(\arabic*)}]

\item $f = \mu_A\xcirc(A \ot\rho)\xcirc (f\ot \mu\ot A)\xcirc (\Delta_{H^{\ot^2}}\ot\eta_A)$,

\item $f$ is a cocycle that satisfies the twisted module condition,

\item $\rho\xcirc (H\ot \eta_A) = \mu_A\xcirc (\rho\ot f) \xcirc (H\ot c\ot H)\xcirc (\Delta\ot\nu)$,

\item $\rho\xcirc (H\ot \eta_A) = \mu_A\xcirc (A\ot f)\xcirc (\nu\ot H)$,

\item $(\mu_A\ot H)\xcirc (A\ot\chi_{\rho})\xcirc(\nu\ot A) = (\mu_A\ot H)\xcirc (A\ot \nu)$,

\end{enumerate}
then

\begin{enumerate}[itemsep=0.7ex,  topsep=1.0ex, label=\emph{(\arabic*)},resume]

\item $\nu$ satisfies condition~1) in Proposition~\ref{condicion 1 sobre 1_E},

\item $\mu_{\mathcal{E}}$ is left and right $A$-linear, associative and normalized with respect to $\nabla_{\!\rho}$,

\item $\nu$ is a preunit of $\mu_{\mathcal{E}}$, $\nabla_{\!\nu}\xcirc \nu = \nu$ and $\nabla_{\!\nu} = \nabla_{\!\rho}$ (consequently, $\imath_{\nu}=\imath_{\rho}$ and $p_{\nu}=p_{\rho}$),

\item $\mu_E$ is left and right $A$-linear, associative and has unit $\eta_E\coloneqq p_{\nu}\xcirc \nu$,

\item the morphism $\jmath'_{\nu}\colon A\to \mathcal{E}$, defined by $\jmath'_{\nu}\coloneqq (\mu_A\ot H)\xcirc (A\ot\nu)$ is left and right $A$-linear, multiplicative and satisfies $\nabla_{\!\nu}\xcirc \jmath'_{\nu}= \jmath'_{\nu}$,

\item the morphism $\jmath_{\nu}\colon A\to E$, defined by $\jmath_{\nu}\coloneqq p_{\nu} \xcirc \jmath'_{\nu}$, is left and right $A$-linear, mul\-ti\-pli\-ca\-ti\-ve and unitary,


\item $\chi_{\rho} = \imath_{\nu}\xcirc \mu_E\xcirc (\gamma \ot \jmath_{\nu})$ and $\mathcal{F}_{\hs f} = \imath_{\nu}\xcirc \mu_E\xcirc (\gamma \ot \gamma)$.

\end{enumerate}

\end{theorem}

\begin{proof} Item~(6) follows from the second equality in~\eqref{cal de gamma}. Thus, by conditions~(1), (3), (4) and~(5), and Propositions~\ref{equivalencia de que F incluido en A times V} and~\ref{cond1preunudad} the tuple $(A,H,\chi_{\rho},\mathcal{F}_{\hs f},\nu)$ is a crossed product system with preunit. Also, by condition~(2) and Proposition~\ref{cocy equiv cocy y..}, the map $\mathcal{F}_{\hs f}$ is a cocycle that satisfies the twisted module condition. So, we can apply Theorem~\ref{prodcruz1} in order to finish the proof.
\end{proof}

From items~(8) and~(9) of Theorem~\ref{weak crossed prod} it follows that
\begin{equation*}\label{gamma 1 = 1}
\gamma\xcirc \eta = p_{\rho}\xcirc (\eta_A\ot \eta) = p_{\nu}\xcirc (\eta_A\ot \eta) = p_{\nu}\xcirc \nu = \eta_E.
\end{equation*}
Moreover, by items~(10) and~(11) of Theorem~\ref{weak crossed prod} and equality~\eqref{cal de gamma'},
\begin{equation}\label{pepe1}
\imath_{\nu} \xcirc \jmath_{\nu} = \nabla_{\!\nu} \xcirc \jmath'_{\nu} = \jmath'_{\nu}\qquad\text{and}\qquad (A\ot \epsilon)\xcirc \jmath'_{\nu} = \mu_A\xcirc (A\ot \rho)\xcirc (A\ot \eta\ot \eta_A).
\end{equation}
So, if $\rho\xcirc (\eta\ot \eta_A) = \eta_A$, then $j_{\nu}$ and $j'_{\nu}$ are monomorphisms. When the hypotheses of this theorem are fulfilled, we say that $E$ is the {\em unitary crossed product of $A$ by $H$ associated with $\rho$ and $f$}.

\smallskip

Given a (not necessarily counitary) right $H$-comodule $\mathcal{B}$ with coaction $\delta_{\mathcal{B}}$, we consider $\mathcal{B}\ot \mathcal{B}$ as a right $H$-comodule via $\delta_{\mathcal{B}\ot \mathcal{B}}\coloneqq (\mathcal{B}\ot \mathcal{B}\ot \mu_H)\xcirc (\mathcal{B}\ot c\ot H) \xcirc (\delta_{\mathcal{B}}\ot \delta_{\mathcal{B}})$.

\begin{proposition}\label{compatibilidad de est de comodule en E con mult}
Under the hypothesis of Theorem~\ref{weak crossed prod} the maps~$\mu_{\mathcal{E}}$ and $\mu_E$ are $H$-colinear.
\end{proposition}

\begin{proof}
For $\mu_{\mathcal{E}}$ this is \cite{FGR}*{Equality (19)}, and for $\mu_E$ it follows easily from the colinearity of~$\mu_{\mathcal{E}}$.
\end{proof}

For each $n\in \mathds{N}_0$, let $\mu_n\colon H^{\ot^n}\to H$ be the map recursively defined by
$$
\mu_0\coloneqq \ide_H\qquad\text{and}\qquad \mu_{n+1}\coloneqq \mu_n\xcirc (\mu\ot H^{\ot^n})\quad\text{for $n\in \mathds{N}_0$.}
$$
We define the maps $u_n\colon H^{\ot^n}\to A$ and $v_n\colon H^{\ot^n}\to A$ by
$$
v_1 = u_1\coloneqq \rho\xcirc (H\ot \eta_A),\quad u_{n+1}\coloneqq u_1\xcirc \mu_n \qquad\text{and}\qquad v_{n+1}\coloneqq u_1\xcirc (H\ot v_n) \quad\text{for $n\in \mathds{N}$,}
$$
where $\rho\colon H\ot A\to A$ denotes the weak measure of $H$ on $A$.

\begin{remark}\label{u_n es idempotente} Since $\rho$ is a weak measure the maps $u_n$ and $v_n$ are idempotent.
\end{remark}

\begin{proposition}\label{alg prop1}
For all $n\in \mathds{N}$, the equality $(v_n\ot\epsilon)\xst v_{n+1} = v_{n+1}$ holds.
\end{proposition}

\begin{proof} For $n = 1$ we have
$$

$$
Clearly the same argumet works for an arbitrary $n$.
\end{proof}

\begin{remark} The equality in Proposition~\ref{equivalencia de que F incluido en A times V}(2) says that $f = f\xst u_2$.
\end{remark}

\begin{remark}\label{epsilon mult * f=f} Let $f\colon H\ot H\to A$ be a map. Arguing as in the proof of Proposition~\ref{equivalencia de que F incluido en A times V} we see that if $u_2\xst f = f$, then $(\epsilon\xcirc \mu\ot f)\xcirc \Delta_{H^{\ot^2}}=f$.
\end{remark}

\begin{remark}\label{a derecha implica aizquierda} If $f\colon H\ot H\to A$ satisfies the twisted module condition, then $f\xst u_2 = v_2\xst f$.
\end{remark}

\begin{definition}\label{aplicacion normal}
 A map $g\colon H\ot H\to A$ is {\em normal} if $g\xcirc (\eta \ot H) = g \xcirc (H\ot \eta) = u_1$.
\end{definition}

\begin{remark}\label{es normal} Assume that $v_2\xst  f = f$. Then,
$$

$$
where the first and third equality hold by definition; the second one, since, by Proposition~\ref{alg prop1} and the hypothesis $(v_1\ot \epsilon)\xst f = (v_1\ot \epsilon)\xst v_2 \xst f = v_2 \xst f = f$; and the last one, by the hypothesis. So, items~(3) and~(4) of Theorem~\ref{weak crossed prod} hold if and only if $f$ is normal.
\end{remark}

In the rest of this subsection we assume that the hypotheses of Theorem~\ref{weak crossed prod} are satisfied.

\begin{lemma}\label{caso particular} The following equalities hold:
\begin{align*}
\mu_{E}\xcirc (\gamma\ot\gamma) &= \mu_{E}\xcirc (\gamma\ot\gamma)\xcirc (\mu\ot\mu)\xcirc (H\ot S\ot H\ot H) \xcirc (H\ot \Delta\xcirc \eta \ot H)\\
&= \mu_{E}\xcirc (\gamma\ot \gamma)\xcirc (\mu\ot \mu)\xcirc (H\ot H\ot S\ot H) \xcirc (H\ot \Delta\xcirc \eta \ot H).
\end{align*}
\end{lemma}

\begin{proof} We prove the first equality and leave the second one to the reader. We have
$$

$$
where the first and last equality hold by Remark~\ref{gama iota y gama gama}; the second one, by Definition~\ref{weak bialgebra}(1) and Propositions~\ref{Delta compuesto con Pi^R multiplicado por algo} and~\ref{Delta eta con S}(2); the third one, by Definition~\ref{weak bialgebra}(3) and Propositions~\ref{Delta compuesto con unidad incluido en H^{hs R}ot H^{hs L}}(1), \ref{fundamental'} and~\ref{caso PiL}; the fourth one, by Remark~\ref{ide ast Pi^R = ide and}; and the fifth one, by Definition~\ref{weak bialgebra}(1).
\end{proof}

\begin{proposition}\label{auxiliar 3} The map $\mu_E\xcirc (E\ot \mu_E)\xcirc (\gamma \ot \jmath_{\nu}\ot \gamma)$ equalize
$$
(H\ot A\ot \mu)\xcirc (\mu \ot c \ot H)\xcirc (H\ot H \ot S\ot A \ot H)\xcirc (H\ot \Delta\xcirc \eta\ot A\ot H)\quad\text{and}\quad \ide_{H\ot A\ot H},
$$
and also equalize
$$
 (H\ot A\ot \mu)\xcirc (\mu \ot c \ot H)\xcirc (H\ot S \ot H\ot A \ot H)\xcirc (H\ot \Delta\xcirc \eta\ot A\ot H)\quad\text{and}\quad \ide_{H\ot A\ot H}.
$$

\end{proposition}

\begin{proof} We prove the first assertion and leave the second one to the reader. We have
$$

$$
where the first and last equality hold by Remark~\ref{gama iota y gama gama}; the second one, by Propositions~\ref{Delta compuesto con unidad incluido en H^{hs R}ot H^{hs L}}(1) and \ref{Delta compuesto con Pi^R multiplicado por algo}; the third one, since $c$ is natural and $\mu_E$ is associative; and fourth one, by Lemma~\ref{caso particular} and the associativity of $\mu_E$.
\end{proof}

\begin{proposition}\label{multiplicacion y gamma} The following assertions hold:

\begin{enumerate}[itemsep=1.0ex, topsep=1.0ex, label=\emph{(\arabic*)}]

\item $\mu_E\xcirc (\gamma \ot \gamma\xcirc \Pi^L)=\gamma \xcirc \mu\xcirc (H\ot \Pi^L)$ and $\mu_E\xcirc (\gamma\xcirc \Pi^L \ot \gamma)=\gamma \xcirc \mu\xcirc (\Pi^L\ot H)$.

\item $\mu_E\xcirc (\gamma \ot \gamma\xcirc \Pi^R)=\gamma \xcirc \mu\xcirc (H\ot \Pi^R)$ and $\mu_E\xcirc (\gamma\xcirc \Pi^R \ot \gamma)=\gamma \xcirc \mu\xcirc (\Pi^R\ot H)$.

\end{enumerate}

\end{proposition}

\begin{proof} (1) \enspace We have
$$

$$
where the first and sixth equality hold because of Proposition~\ref{gama iota y gama gama}; the second one, by Propositions~\ref{Delta compuesto con unidad incluido en H^{hs R}ot H^{hs L}}(1), \ref{Delta compuesto con Pi^R}(3), \ref{fundamental'} and~\ref{caso PiL}; the third one, since $f$ is normal by Remark~\ref{es normal}; the fourth one, by Definition~\ref{weak bialgebra}(1); the fifth one, by the definition of $u_1$; and the last one, by Theorem~\ref{weak crossed prod}(11). So the first equality is true. A similar argument proves the second one.

\smallskip

\noindent (2)\enspace Mimic the proof of item~(1).
\end{proof}

\begin{proof} This follows from Propositions~\ref{conmutatividad debil} and~\ref{multiplicacion y gamma}.
\end{proof}

\begin{corollary}\label{gamma xst (gamma xcirc Pi^R) = gamma}
The equalities $(\gamma \xcirc \Pi^{\hs L}) \xst \gamma  = \gamma \xst (\gamma \xcirc \Pi^{\hs R}) = \gamma $ hold.
\end{corollary}

\begin{proof} By Proposition~\ref{multiplicacion y gamma}(1) and Remark~\ref{ide ast Pi^R = ide and} we have $\gamma  = \gamma \xcirc (\Pi^{\hs L} \xst \ide ) = (\gamma \xcirc \Pi^{\hs L})\xst \gamma$. A similar  argument shows that $\gamma = \gamma \xst (\gamma \xcirc \Pi^{\hs R})$.
\end{proof}

\subsection[Right $H$-comodule algebras]{Right $\bm{H}$-comodule algebras}\label{Right comodule algebras}

\begin{proposition}\label{wbialgebras} Let $\mathcal{B}$ be a right $H$-comodule which is also an unitary algebra. If $\mu_{\mathcal{B}}$ is an $H$-colinear map, then the following assertions are equivalent:

\begin{enumerate}[itemsep=0.7ex,  topsep=1.0ex, label=\emph{(\arabic*)}]

\item $(\mathcal{B}\ot\Delta)\xcirc \delta_{\mathcal{B}}\xcirc \eta_{\mathcal{B}} = (\mathcal{B}\ot \mu\ot H)\xcirc (\delta_{\mathcal{B}}\ot \Delta)\xcirc (\eta_{\mathcal{B}}\ot \eta)$.

\item $(\mathcal{B}\ot\Delta)\xcirc \delta_{\mathcal{B}}\xcirc \eta_{\mathcal{B}}  = (\mathcal{B}\ot \mu\ot H)\xcirc (\mathcal{B}\ot c\ot H)\xcirc (\delta_{\mathcal{B}}\ot \Delta)\xcirc (\eta_{\mathcal{B}}\ot \eta)$.

\item $(\mathcal{B}\ot \ov{\Pi}^R)\xcirc \delta_{\mathcal{B}} = (\mu_{\mathcal{B}}\ot H)\xcirc (\mathcal{B}\ot \delta_{\mathcal{B}})\xcirc (\mathcal{B}\ot \eta_{\mathcal{B}})$.

\item $(\mathcal{B}\ot \Pi^{\hs L})\xcirc \delta_{\mathcal{B}} = (\mu_{\mathcal{B}}\ot H)\xcirc (\mathcal{B}\ot c)\xcirc (\delta_{\mathcal{B}}\ot \mathcal{B})\xcirc (\eta_{\mathcal{B}}\ot \mathcal{B})$.

\item $(\mathcal{B}\ot \ov{\Pi}^R)\xcirc \delta_{\mathcal{B}} \xcirc \eta_{\mathcal{B}}=\delta_{\mathcal{B}}\xcirc \eta_{\mathcal{B}}$.

\item $(\mathcal{B}\ot \Pi^{\hs L})\xcirc \delta_{\mathcal{B}} \xcirc \eta_{\mathcal{B}}=\delta_{\mathcal{B}}\xcirc \eta_{\mathcal{B}}$.

\end{enumerate}

\end{proposition}

\begin{proof} See~\cite{CG}.
\end{proof}

\begin{definition}\label{def: comodulo  algebra} An unitary algebra $\mathcal{B}$, which is also a right $H$-comodule, is a {\em right $H$-comodule algebra} if $\mu_{\mathcal{B}}$ is $H$-colinear and the equivalent statements of the previous proposition are satisfied.
\end{definition}

\smallskip

The following result and its proof were kindly communicated to us by Jos\'e Nicanor Alonso \'Alvarez y Ram\'on Gonz\'alez Rodr\'iguez. It is is interesting in itself and allows to simplify the proof of Proposition~\ref{E es H-comodulo algebra debil}, but only under the assumption that $H$ is a weak Hopf algebra.

\begin{proposition}\label{suf para comodulo algebra} Let $\mathcal{B}$ be a right $H$-comodule which is also an unitary algebra. Assume that $\mu_{\mathcal{B}}$ is an $H$-colinear map. If $H$ is a weak Hopf algebra, then the equivalent items of the previous proposition are fulfilled.
\end{proposition}

\begin{proof}
By the hypotheses, Proposition~\ref{mu Pi^R, etc}(1) and Definition~\ref{weak Hopf Algebra}(1), we have
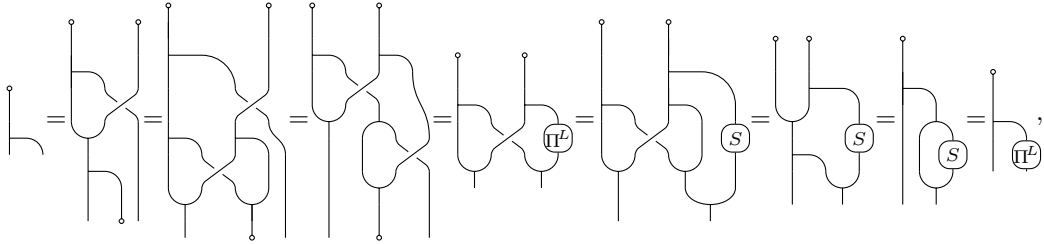
\begin{equation*}

\end{equation*}
as desired.
\end{proof}

\begin{proposition}[Comodule algebra structure on an unitary crossed product $E$]\label{E es H-comodulo algebra debil} Each unitary crossed product $E$ is a weak $H$-comodule algebra via the coaction introduced in Remark~\ref{coaccion en A ot H}.
\end{proposition}

\begin{proof} By Proposition~\ref{compatibilidad de est de comodule en E con mult} we know that $\mu_E$ is $H$-colinear. We next prove that the equality in item~(1) of Proposition~\ref{wbialgebras} is satisfied. That is
\begin{equation}\label{eq: compatibilidad de est de comodule en E con mult}
(E\ot \Delta)\xcirc \delta_E\xcirc \eta_E = (E\ot\mu\ot H)\xcirc(\delta_E\ot \Delta)\xcirc (\eta_E \ot\eta).
\end{equation}
By the fact that $\nabla_{\!\rho}\xcirc \nu = \nu$, condition~(2) in Proposition~\ref{condicion 1 sobre 1_E} and the first equality in Definition~\ref{weak bialgebra}(2), we have
$$

$$
This finishes the proof, since the left side in equality~\eqref{eq: compatibilidad de est de comodule en E con mult} is the map represented by the first diagram, and the right side is the map represented by the last diagram.
\end{proof}

\begin{definition}\label{integral} Let $\mathcal{B}$ be a right $H$-comodule algebra. An {\em integral} of $\mathcal{B}$ is a morphism of right $H$-comodules $\gamma\colon H\to \mathcal{B}$. If moreover $\gamma(1) = 1_{\mathcal{B}}$, then $\gamma$ is a {\em total integral}.
\end{definition}

\begin{definition}\label{integral conv inv} An integral $\gamma$ is {\em convolution invertible} if there exists a morphism $\gamma^{-1}\colon H\to \mathcal{B}$ such that
$$
\gamma^{-1}\xst \gamma = \gamma\xcirc \Pi^{\hs R},\quad \gamma\xst \gamma^{-1} = \gamma\xcirc \Pi^L\quad\text{and}\quad (\gamma\xcirc \Pi^{\hs R})\xst \gamma^{-1} = \gamma^{-1}.
$$
Clearly $\gamma^{-1}$ is unique.
\end{definition}

\begin{definition}\label{cleft} Let $\mathcal{B}$ be a right $H$-comodule algebra and $\jmath\colon A\to \mathcal{B}$ an algebra monomorphism. We say that $(\mathcal{B},\jmath)$ is an {\em extension} of $A$ by $H$ if $\jmath$ is the equalizer of $\delta_{\mathcal{B}}$ and $(\mathcal{B}\ot \Pi^L)\xcirc \delta_{\mathcal{B}}$, and we say that a extension $(\mathcal{B},\jmath)$ is {\em $H$-cleft} if there exists a convolution invertible total integral $\gamma\colon H\to \mathcal{B}$ such that $\gamma\xcirc \Pi^L$ factorizes through $\jmath$.
\end{definition}

In the sequel $\gamma$ will be called a {\em cleaving map} associated with the $H$-cleft extension $(\mathcal{B},\jmath)$.

\smallskip

\begin{definition}\label{morfismo de extensiones}
Let $(\mathcal{B},\jmath)$ and $(\mathcal{B}',\jmath')$ be extensions of $A$ by $H$. An arrow $\Phi\colon \mathcal{B}\to \mathcal{B}'$ is a {\em morphism of extensions} if $\Phi$ is a $H$-colinear algebra morphism and $\Phi\xcirc \jmath = \jmath'$.
\end{definition}

Let $\Phi$ be as in the previous definition. If $(\mathcal{B},\jmath)$ is an $H$-cleft extension with cleaving map $\gamma$ and inverse $\gamma^{-1}$, then $(\mathcal{B}',\jmath')$ is an $H$-cleft extension with cleaving map $\Phi\xcirc \gamma$ and inverse $\Phi\xcirc \gamma^{-1}$. Moreover, in this case, $\Phi$ is an isomorphism.

\begin{definition}\label{equivalencia de extensiones}
Two extensions of $A$ by $H$ are {\em equivalent} if they are isomorphic.
\end{definition}

\section{Equivalence of weak crossed products}\label{Sec: Equivalence of weak crossed products}

Let $H$ be a weak bialgebra and let $A$ be a unitary algebra. Let $\rho$ and $\rho'$ be weak measures of $H$ on $A$, and let $f\colon H\ot H\to A$ and $f'\colon H\ot H\to A$ be maps. Assume that both pairs $(\rho,f)$ and $(\rho',f')$ satisfy the hypotheses of Theorem~\ref{weak crossed prod}. In this section we set $\mathcal{E}\coloneqq A\ot_{\rho}^f H$, $E\coloneqq A\times_{\rho}^f H$, $\mathcal{E}'\coloneqq A\ot_{\rho'}^{f'} H$ and $E'\coloneqq A\times_{\rho'}^{f'} H$. Recall that the preunits of $\mathcal{E}$ and $\mathcal{E}'$ are the maps defined by $\nu = \chi_{\rho}\xcirc (\eta\ot \eta_A)$ and $\nu' = \chi_{\rho'}\xcirc (\eta\ot \eta_A)$, respectively, while the units of $E$ and $E'$ are the maps $\eta_E = p_{\rho}\circ \nu$ and $\eta_{E'} = p_{\rho'}\circ \nu'$. In the sequel we let $u_1\colon H\to A$ and $u'_1\colon H\to A$  denote the maps defined by $u_1\coloneqq \rho\xcirc (H\ot \eta_A)$ and $u'_1\coloneqq \rho'\xcirc (H\ot \eta_A)$. Moreover, for each map $\phi\colon H\to A$, we let $L(\phi)\colon A\ot H\to A\ot H$ denote the left $A$-linear and right $C$-colinear map defined by $L(\phi)\coloneqq (\mu_A\ot H)\xcirc (A\ot \phi\ot H)\xcirc (A\ot \Delta)$.

\begin{theorem}\label{caracterizacion de prod cruz equiv} If $\phi\colon H\to A$ is a map that satisfies
\begin{enumerate}[itemsep=0.7ex,  topsep=1.0ex, label=\emph{(\arabic*)}]

\item $\phi= u_1\xst \phi = \phi\xst u'_1$,

\item There exists $\phi'\colon H\to A$ such that $\phi \xst \phi' = u_1$ and $\phi'\ast \phi =u'_1$,

\item $(\phi \ot H)\xcirc \Delta \xcirc \eta = \nu'$,

\item $\mu_A\xcirc (A\ot \phi) \xcirc \chi_{\rho} =\mu_A\xcirc (\phi\ot \rho') \xcirc (\Delta\ot A)$,

\item $\mu_A\xcirc (A\ot \phi)\xcirc \mathcal{F}_{\hs f} = \mu_A\xcirc (\mu_A\ot f')\xcirc (A\ot \chi_{\rho'}\ot H)\xcirc (\phi\ot H \ot \phi\ot H) \xcirc (\Delta\ot \Delta)$,

\end{enumerate}
then the map $\Phi_{\!_{\phi}}\colon E\to E'$, defined by $\Phi_{\!_{\phi}}\coloneqq p_{\rho'}\xcirc L(\phi) \xcirc \imath_{\rho}$, is a left $A$-linear and right $H$-co\-lin\-ear isomorphism of unitary algebras. Conversely, if $\Phi\colon E\to E'$ be a left $A$-linear and right $H$-colinear isomorphism of unitary algebras, then the map $\phi_{\!_{\Phi}}\coloneqq (A\ot \epsilon)\xcirc \imath_{\rho'}\xcirc \Phi \xcirc p_{\rho}\xcirc (\eta_A\ot H)$ satisfy statements~(1)--(5). Moreover, the correspondences $\Phi\mapsto \phi_{\!_{\Phi}}$ and  $\phi\mapsto \Phi_{\!_{\phi}}$ are inverse one of each other.
\end{theorem}

In order to prove this result, we first establish a sequence of Lemmas.

\begin{lemma}\label{composicion con nabla} For each map  $\phi\colon H\to A$, the following facts hold:

\begin{enumerate}[itemsep=0.7ex,  topsep=1.0ex, label=\emph{(\arabic*)}]

\item $\phi= u_1\ast \phi$ if and only if $L(\phi)=L(\phi)\xcirc \nabla_{\!\rho}$.

\item $\phi= \phi \ast u_1'$ if and only if $L(\phi)= \nabla_{\!\rho'}\xcirc L(\phi)$.

\end{enumerate}

\end{lemma}

\begin{proof} Left to the reader.  \end{proof}

\begin{lemma}\label{composicion de phi y phi'}
Let $\phi\colon H\to A$ and $\phi'\colon H\to A$ be maps and let $\Phi\colon E\to E'$ and $\Phi'\colon E'\to E$ be the maps  $\Phi\coloneqq \Phi_{\!_{\phi}}$ and $\Phi'\coloneqq \Phi_{\!_{\phi'}}$. Assume that $\phi=u_1\ast \phi=\phi\ast u'_1$ and $\phi'=u'_1\ast \phi'=\phi'\ast u_1$.
Then the following facts hold:

\begin{enumerate}[itemsep=0.7ex,  topsep=1.0ex, label=\emph{(\arabic*)}]

\item $\phi\ast \phi'=u_1$ if and only if $\Phi'\xcirc \Phi=\ide_E$.

\item $(\phi\ot H)\circ \Delta \circ \eta=\nu'$ if and only if $\Phi\circ \eta_E=\eta_{E'}$.

\item  The maps $\chi\coloneqq \chi_{\rho}$ and $\chi'\coloneqq \chi_{\rho'}$ have the following property:
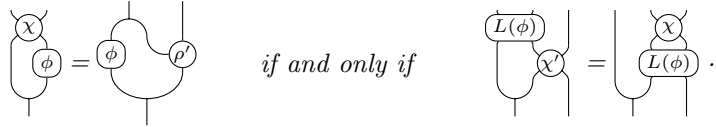
\begin{equation*}
\qquad

\end{equation*}

\item If the equivalent conditions in item~(3) are fulfilled, then
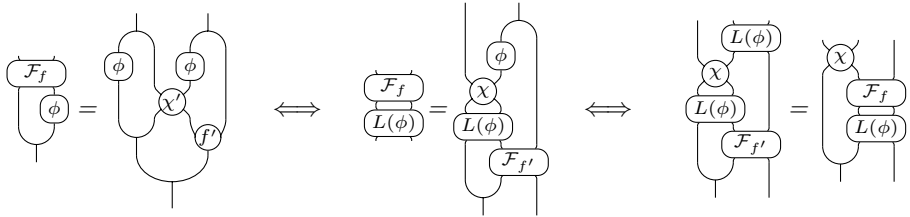
\begin{equation*}
\qquad%
\end{equation*}

\end{enumerate}

\end{lemma}

\begin{proof}
We prove that last item and left the other ones, which are easier, to the reader. By item~(3), the definition of $L(\phi)$ and Remark~\ref{aux}(2), we have
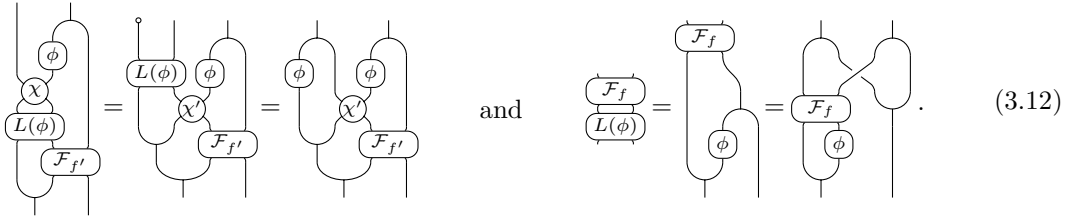
\begin{equation}\label{us1}
%
\end{equation}
Using this, the coassociativity of $\Delta$, the fact that~$c$ is natural and Remark~\ref{aux}(1), we obtain that, if the first equality in item~(4) of the statement is true, then
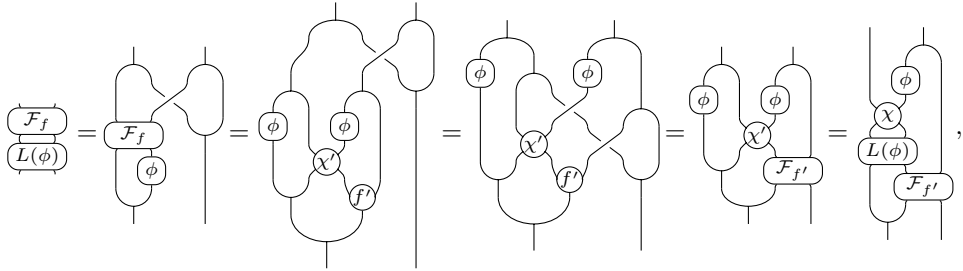
\begin{equation*}
%
\end{equation*}
which proves the second equality in item~(4) of the statement. Conversely, if this equality is true, then, again by~\eqref{us1}, we have
\begin{equation*}
%
\end{equation*}
which proves the first equality in item~(4) of the statement. Assume again that the second equal\-i\-ty in the state\-ment is true. Then, we have
\begin{equation*}
%
\end{equation*}
which proves that the last equality in item~(4) of the statement is true. Finally, if this equality holds, then
\begin{equation*}
%
\end{equation*}
which proves the second equality in item~(4) of the statement.
\end{proof}

\begin{proof}[Proof of Theorem~\ref{caracterizacion de prod cruz equiv}] Assume that $\phi$ satisfies statements~1)--5) and set $\Phi\coloneqq \Phi_{\!_{\phi}}$. The map~$\Phi$ is left $A$-linear and right $H$-colinear because $p_{\rho'}$, $L(\phi)$ and $\imath_{\rho}$ are left $A$-linear and right $H$-colinear maps and
$$
\Phi = p_{\rho'}\xcirc L(\phi) \xcirc \imath_{\rho}.
$$
We next prove that it is also an isomorphism of unitary algebras. By Lemma~\ref{composicion con nabla} and the fact that $\mu_{\mathcal{E}}$ is normalized with respect to $\nabla_{\!\rho}$ and $\mu_{\mathcal{E'}}$ is normalized with respect to $\nabla_{\!\rho'}$, we have
\begin{equation*}

\end{equation*}
Consequently $\Phi$ is multiplicative if and only if $L(\phi)$ is. But this is true, because, by items~(3) and~(4) of Lemma~\ref{composicion de phi y phi'} and the fact that $L(\phi)$ is left $A$-linear,
\begin{equation*}
%
\end{equation*}
where $\chi\coloneqq \chi_{\rho}$ and $\chi'\coloneqq \chi_{\rho'}$. It remains to prove that the morphism $\Phi$ is unitary and invertible. But $\Phi$ is unitary by item~(2) of Lemma~\ref{composicion de phi y phi'}, and using item~(1) of Lemma~\ref{composicion de phi y phi'} we obtain that $\Phi$ is invertible (with inverse given by the map $\Phi_{\!_{\phi'}}\coloneqq p_{\rho}\xcirc L(\phi')\xcirc \imath_{\rho'}$).

\smallskip

Conversely, assume that $\Phi\colon E \to E'$ is an unitary algebra isomorphism that is left $A$-linear and right $H$-colinear and set $\phi = \phi_{\!_{\Phi}}\coloneqq (A\ot \epsilon)\xcirc \imath_{\rho'}\xcirc \Phi \xcirc p_{\rho}\xcirc (\eta_A\ot H)$. Since $ p_{\rho}\xcirc\nabla_{\!\rho}=p_{\rho}$ and $p_{\rho}$, $\Phi$ and $i_{\rho'}$ are left $A$-linear and left $H$-colinear,
\begin{equation*}

\end{equation*}
which by Lemma~\ref{composicion con nabla} proves that $\phi= u_1\ast \phi$. A similar computation proves that $\phi= \phi\ast u'_1$. So, item~(1) is fulfilled. Items~(2) and~(3) follow from items~(1) and~(2) of Lemma~\ref{composicion de phi y phi'}, since
\begin{equation}\label{us3}
%
\end{equation}
where $\Phi'$ is the inverse of $\Phi$ and $\phi'\coloneqq \phi_{\!_{\Phi'}}\coloneqq (A\ot \epsilon)\xcirc \imath_{\rho}\xcirc \Phi' \xcirc p_{\rho'}\xcirc (\eta_A\ot H)$. Note now that
$$
L(\phi)\xcirc \nu = L(\phi)\xcirc \nabla_{\!\rho}\xcirc (\eta_A\ot \eta) = L(\phi)\xcirc (\eta_A\ot \eta) = (\phi\ot H)\xcirc \Delta\xcirc \eta = \chi_{\rho'}\xcirc (\eta\ot\eta_A) = \nu',
$$
where the second equality holds by Lemma~\ref{composicion con nabla}; the third one, by the definition of $L(\phi)$; the fourth one, by  Lemma~\ref{composicion de phi y phi'}(2); and the last one, by the first equality in~\eqref{cal de gamma}. Hence,
$$
L(\phi)\xcirc \jmath'_{\nu}\! =\! L(\phi) \xcirc (\mu_A\ot H) \xcirc (A\ot\nu)\! =\! (\mu_A\ot H)\xcirc (A\ot L(\phi)) \xcirc (A\ot \nu)\! =\! (\mu_A\ot H) \xcirc (A\ot \nu')\! =\! \jmath'_{\nu'},
$$
where the second equality holds because $L(\phi)$ is left $A$-linear. Consequently, by the first equality in Theorem~\ref{prodcruz1}(8) and the fact that $L(\phi)$ is multiplicative and both $L(\phi)$ and $\mu_{\mathcal{E}}$ are left $A$-linear, we have
\begin{equation}\label{us2}

\end{equation}
which by Lemma~\ref{composicion de phi y phi'}(3) implies that item~(4) is fulfilled. Using now equality~\eqref{us2} and the fact that $L(\phi)$ is multiplicative, we obtain that
\begin{equation*}
%
\end{equation*}
which by Lemma~\ref{composicion de phi y phi'}(4) proves that item~(5) is true.

\smallskip

Finally, the first calculations made out in~\eqref{us3} and the fact that
$$
\phi = (A\ot \epsilon)\xcirc L(\phi)\xcirc (\eta_A\ot H) = (A\ot \epsilon)\xcirc \imath_{\rho'} \xcirc \Phi_{\phi} \xcirc p_{\rho} \xcirc (\eta_A\ot H) = \phi_{\Phi_{\phi}},
$$
%
%
%
show that the correspondences $\Phi\mapsto \phi_{\!_{\Phi}}$ and  $\phi\mapsto \Phi_{\!_{\phi}}$ are inverses one of each other.
\end{proof}

\begin{definition}\label{productos debiles cruzados equivalentes} Two crossed products $E\coloneqq A\times_{\rho}^{f} H$ and $E'\coloneqq A\times_{\rho'}^{f'} H$ are said to be {\em equivalent} if the exists an algebra isomorphism $\Phi\colon E\to E'$ that is left $A$-linear and right $H$-colinear.
\end{definition}

\begin{remark}\label{la equivalencia extiende la del caso coconmutativo} Theorem~\ref{caracterizacion de prod cruz equiv} shows that the notion of equivalence of crossed product of algebras by weak Hopf algebras given in this section reduce to the one introduced in \cite{AVR} above Proposition~3.2.
\end{remark}

\section{Weak crossed products of weak module algebras}\label{subsection: Weak crossed products of weak module algebras by weak bialgebras in which the unit cocommutes}
Let $H$ be a weak bialgebra, $A$ and algebra and $\rho\colon H\ot A\to A$ a map. In this section we study the weak crossed products of $A$ with $H$ in which $A$ is a left weak $H$-module algebra.

\begin{proposition}\label{wbialgebras'} Let $A$ be an unitary algebra. If $\rho$ satisfies

\begin{enumerate}[itemsep=0.7ex,  topsep=1.0ex, label=\emph{(\arabic*)}]

\item $\rho\xcirc (\eta\ot A) = \ide_A$,

\item $\rho \xcirc (H\ot \mu_A) = \mu_A\xcirc (\rho\ot \rho)\xcirc (H\ot c\ot A)\xcirc (\Delta\ot A\ot A)$,

\item $\rho\xcirc (\mu\ot \eta_A) = \rho \xcirc (H\ot \rho)\xcirc (H\ot H\ot\eta_A)$,

\end{enumerate}
then the following assertions are equivalent:

\begin{enumerate}[itemsep=0.7ex,  topsep=1.0ex, label=\emph{(\arabic*)},resume]

\item $\rho \xcirc (\Pi^{\hs L}\ot A) = \mu\xcirc (\rho\ot A)\xcirc (H\ot \eta_A\ot A)$.

\item $\rho \xcirc (\ov{\Pi}^L\ot A) = \mu\xcirc c\xcirc (\rho\ot A)\xcirc (H\ot \eta_A\ot A)$.

\item $\rho \xcirc (\Pi^{\hs L}\ot \eta_A) = \rho\xcirc (H\ot \eta_A)$.

\item $\rho \xcirc (\ov{\Pi}^L\ot \eta_A) = \rho\xcirc (H\ot \eta_A)$.

\item $\rho\xcirc (H\ot \rho) \xcirc (H^{\ot^2}\ot \eta_A) = (\rho\ot \epsilon)\xcirc (H\ot c)\xcirc (H\ot \mu\ot A)\xcirc (\Delta\ot H\ot \eta_A)$.

\item $\rho\xcirc (H\ot \rho) \xcirc (H^{\ot^2}\ot \eta_A) = (\epsilon\ot A)\xcirc (\mu\ot \rho)\xcirc (H\ot c\ot A)\xcirc (\Delta\ot H\ot \eta_A)$.

\end{enumerate}

\end{proposition}

\begin{proof}
See~\cite{CG}.
\end{proof}

When $\rho$ satisfies items~(1), (2), (3) and~(4) of Proposition~\ref{wbialgebras'} we say that $A$ is a {\em left weak $H$-module algebra} and $\rho$ is a {\em weak left action of $H$ on $A$.} If also $\rho\xcirc (\mu\ot A) = \rho \xcirc (H\ot \rho)$, then we say that $\rho$ is an {\em action} and $A$ is a {\em left $H$-module algebra}.

\smallskip

Like Proposition~\ref{suf para comodulo algebra}, the following result was communicated to us by Jos\'e Nicanor Alonso \'Alvarez y Ram\'on Gonz\'alez Rodr\'iguez. It is is interesting in itself and allows to simplify the proof of Proposition~\ref{es un H modulo algebra debil}, but only under the assumption that $H$ is a weak Hopf algebra.

\begin{proposition} If $H$ is a weak Hopf algebra and $\rho$ satisfies the hypotheses of the above proposition, then the equivalent conditions~4)--9) of that proposition are satisfied.
\end{proposition}

\begin{proof}
By the hypotheses, Proposition~\ref{delta Pi^R, etc} and Definition~\ref{weak Hopf Algebra}(1), we have
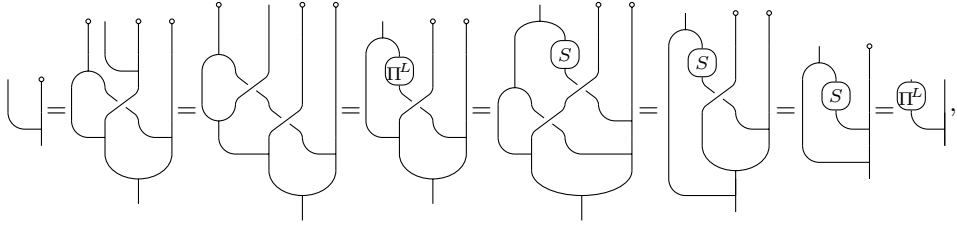
\begin{equation*}

\end{equation*}
as desired.
\end{proof}

\begin{remark}\label{accion de PiL y comultipliacion} By the definition of $\ov{\Pi}^{\hs L}$, the fact that the fact that $c$ is natural and $\Delta$ is counitary, and Propositions~\ref{Delta compuesto con unidad incluido en H^{hs R}ot H^{hs L}}(1) and~\ref{Delta compuesto con Pi^R multiplicado por algo},
$$
%
$$
\end{remark}

\begin{proposition} We have $\nabla_{\nu} = (\rho\ot\mu) \xcirc (H\ot c\ot H)\xcirc (\Delta\ot A\ot H)\xcirc (\eta\ot A\ot H)$.
\end{proposition}

\begin{proof} By Remark~\ref{accion de PiL y comultipliacion}, Theorem~\ref{weak crossed prod}(8), equality~\eqref{calculo de nabla rho}, the fact that $c$ is natural, the equalities in items~(5) and~(7) of Proposition~\ref{wbialgebras'}, and Proposition~\ref{Delta compuesto con unidad incluido en H^{hs R}ot H^{hs L}}(2), we have
$$

$$
as desired.
\end{proof}

From here to the end of this subsection $A$ is a left weak $H$-module algebra and $E$ is the unitary crossed product of $A$ by $H$ associated with $\rho$ and a map $f\colon H\ot H\to A$. Thus, we assume that the hypotheses of Theorem~\ref{weak crossed prod} are fulfilled. In particular $\nu$, $\jmath_{\nu}$ an $\gamma$ are as in that theorem. In the sequel we study the consequences that $A$ is a left weak $H$-module algebra. By the equality in Proposition~\ref{wbialgebras'}(3), we have $v_n = u_n$ for all $n\in \mathds{N}$.

\begin{proposition}\label{el egalizador es A}
The pair $(A,\jmath_{\nu})$ is the equalizer of $\delta_E$ and $(E\ot \Pi^L)\xcirc \delta_E$.
\end{proposition}

\begin{proof}
Let $g\colon X \to E$ be an arrow such that $\delta_E\xcirc g =(E\ot \Pi^L)\xcirc\delta_E\xcirc g$. Then,
$$

$$
where the first equality holds since $\nabla_{\!\nu}$ is $H$-colinear, $\Delta$ is counitary and $\nabla_{\!\nu}\xcirc \imath_{\nu}=\imath_{\nu}$; the second one, since $\nabla_{\!\nu}=\nabla_{\!\rho}$ and $\delta_E\xcirc g = (E\ot \Pi^L) \xcirc \delta_E\xcirc g$; the third one, by the very definition of $\Pi^L$; the fourth one, by the fact that $\Delta$ is counitary and Propositions~\ref{Delta compuesto con unidad incluido en H^{hs R}ot H^{hs L}}(1) and~\ref{Delta compuesto con Pi^R multiplicado por algo}; the fifth one, by Proposition~\ref{Delta compuesto con unidad incluido en H^{hs R}ot H^{hs L}}(2) and the equalities in items~(3) and~(5) of Proposition~\ref{wbialgebras'}; and the sixth one by the associativity of $\mu_A$ and the definition of $\jmath'_{\nu}$. Since $\jmath'_{\nu}=\imath_{\nu}\xcirc \jmath_{\nu}$ and $\imath_{\nu}$ is a monomorphism, this implies that
$$
g=  \jmath_{\nu} \xcirc \mu_A \xcirc (A\ot \rho) \xcirc (\imath_{\nu} \ot A) \xcirc (g\ot \eta_A).
$$
So, $g$ factorizes throught $\jmath_{\nu}$. Since $\jmath_{\nu}$ is a monomorphism because $(A\ot \epsilon)\xcirc \imath_{\nu}\xcirc \jmath_{\nu}= \ide_A$, this proves the assertion.
\end{proof}

\begin{proposition}\label{auxiliar4} The equality $\gamma \xcirc \Pi^{\hs L}=\jmath_{\nu}\xcirc \rho\xcirc (H\ot \eta_A)$ holds.
\end{proposition}

\begin{proof}
To begin note that by items~(10) and~(11) of Theorem~\ref{weak crossed prod}, we have $\imath_{\nu}\xcirc \jmath_{\nu}=\nabla_{\nu}\xcirc \jmath'_{\nu}=\jmath'_{\nu}$. Using this, Theorem~\ref{weak crossed prod}(8), the first and last equality in~\eqref{cal de gamma}, Proposition~\ref{Delta compuesto con Pi^R multiplicado por algo}, the equalities in items~(3) and~(4) of Proposition~\ref{wbialgebras'}, and the definition of $\jmath'_{\nu}$, we obtain that
$$

$$
Since $\imath_{\nu}$ is a monomorphism, the statement follows from these equalities.
\end{proof}

\begin{proposition}\label{conmutatividad de gamma e i} The equality $\mu_E\xcirc (\gamma\ot \jmath_{\nu})\xcirc (\Pi^{\hs R}\ot A) = \mu_E\xcirc (\jmath_{\nu}\ot \gamma)\xcirc (A\ot \Pi^{\hs R})\xcirc c$ holds.
\end{proposition}

\begin{proof} We have
$$
%
$$
where the first and fifth equalities hold by Remark~\ref{gama iota y gama gama}; the second one, by Propositions~\ref{algunas composiciones} and~\ref{mu delta Pi^R, etc}; the third one, by the equality in Proposition~\ref{wbialgebras'}(5); the fourth one, by Theorem~\ref{weak crossed prod}(11) and the associativity of $\mu_E$; and the last one, again by Theorem~\ref{weak crossed prod}(11).
\end{proof}

\section[Weak crossed products with invertible cocycle are cleft]{Weak crossed products with invertible cocycle are cleft}\label{cociclo inversible}

In this section we assume that $H$ is a weak Hopf algebra, $A$ is a left weak $H$-module algebra and $E$ is the unitary crossed product of $A$ by $H$ associated with a weak left action $\rho$ and a map $f\colon H\ot H\to A$. Thus, we assume that the hypotheses of Theorem~\ref{weak crossed prod} and the equalities in the items of Proposition~\ref{wbialgebras'} are fulfilled. The aim of this section is to prove that the crossed products with invertible cocycle are $H$-cleft extensions.

\subsection{Regular maps}

\begin{proposition}\label{g*u_n} Let $g,g'\in \Hom(H^{\ot^n},A)$. If $g\xst g' = g'\xst g = u_n$, then

\begin{enumerate}[itemsep=0.7ex,  topsep=1.0ex, label=\emph{(\arabic*)}]

\item $g\xst u_n = u_n\xst g$ and $g'\xst u_n  = u_n\xst g'$,

\item $(u_n\xst g)\xst (u_n\xst g') = u_n$ and $(u_n\xst g')\xst (u_n\xst g) = u_n$.

\end{enumerate}

\end{proposition}

\begin{proof} By symmetry we only must prove the first equalities en items~(1) and~(2). But
$$
g\xst u_n = g\xst (g'\xst g) = (g\xst g')\xst g = u_n\xst g\quad\text{and}\quad (u_n\xst g)\xst (u_n\xst g') = (u_n\xst u_n)\xst (g\xst g') = u_n,
$$
as desired.
\end{proof}

\begin{proposition}\label{alg prop} Let $F_1$, $F_2$ and $F_{\!\rho}$ be the morphisms from $H^{\ot^3}$ to $A$ defined by $F_1\coloneqq f\xcirc (\mu\ot H)$, $F_2\coloneqq f\xcirc (H\ot \mu)$ and $F_{\!\rho}\coloneqq \rho\xcirc (H\ot f)$. The following facts hold:
\begin{enumerate}[itemsep=0.7ex,  topsep=1.0ex, label=\emph{(\arabic*)}]

\item $F_1 \xst  u_3= F_1$ and $u_3 \xst  F_1 = F_1$,

\item $F_2 \xst  u_3= F_2$ and $u_3 \xst  F_2 = F_2$,

\item $F_{\!\rho} \xst  u_3= F_{\!\rho}$ and $u_3 \xst  F_{\!\rho} = F_{\!\rho}$.

\end{enumerate}

\end{proposition}

\begin{proof} (1)\enspace By the fact that $f * u_2 = f$ and Definition~\ref{weak bialgebra}(1),
$$

$$
Similarly, $u_3 \xst  F_1 = F_1$.

\smallskip

\noindent (2)\enspace Mimic the proof of item~(1).

\smallskip

\noindent (3)\enspace By the fact that $f * u_2 = f$ and the equalities in item~(2) and~(3) of Proposition~\ref{wbialgebras'},
$$
%
$$
Similarly, $u_3 \xst  F_{\!\rho} = F_{\!\rho}$.
\end{proof}

\begin{proposition}\label{vale} The equality $f\xcirc (H\ot \mu\xcirc (H\ot \Pi^L)) = f\xcirc (H\ot \mu\xcirc (H\ot \ov{\Pi}^L))$ holds.
\end{proposition}

\begin{proof} It suffices to prove that
\begin{equation}\label{dos igualdades}
f\xcirc (H\ot \mu)\xcirc (H\ot H\ot \Pi^L) =  \mu_A\xcirc (A\ot u_1)\xcirc (\mathcal{F}_f\ot H) =  f\xcirc (H\ot \mu)\xcirc (H\ot H\ot \ov{\Pi}^L).
\end{equation}
The first equality is true since
$$
%
$$
where the first equality holds because $f = v_2*f$ by Remark~\ref{a derecha implica aizquierda}; the second one, by the fact that $f$ is normal and the equality in Definition~\ref{weak bialgebra}(1); the third one, by Propositions~\ref{Delta compuesto con unidad incluido en H^{hs R}ot H^{hs L}}(1), \ref{Delta compuesto con Pi^R}(3), \ref{fundamental'} and~\ref{caso PiL}; the fourth one, by the cocycle condition; The fifth one by Proposition~\ref{caso PiL} and the fact that $f$ is normal; and the last one, by the fact that $f$ is normal and the equalities in items~(3) and~(6) of Proposition~\ref{wbialgebras'}. The proof of the second equality in~\eqref{dos igualdades} is similar, but we must use that $\ov{\Pi}^L = \Pi^R\xcirc \ov{\Pi}^L$ and Propositions~\ref{mu delta Pi^R, etc} and~\ref{fundamental'} instead of Propositions~\ref{Delta compuesto con unidad incluido en H^{hs R}ot H^{hs L}}(1), \ref{Delta compuesto con Pi^R}(3), \ref{fundamental'} and~\ref{caso PiL}.
\end{proof}

\begin{remark}\label{rem1} We have
$$

$$
where the first equality holds by items~(3), (6) and~(9) of Proposition~\ref{wbialgebras'}; the second one, by condition~(1) of Definition~\ref{weak bialgebra}; the third one, by the coassociativity of $\Delta_{H\ot H}$; and the last one, since $c$ is natural.
\end{remark}

\begin{proposition}\label{F_epsilon*u_3 = u_3*F_epsilon} Let $F_{\!\epsilon}\colon H^{\ot^3}\to A$ be the morphism defined by $F_{\!\epsilon}\coloneqq f\ot \epsilon$. Then:
$$
F_{\!\epsilon}\xst u_3 = u_3\xst F_{\!\epsilon} = f\xcirc (H\ot \mu)\xcirc (H^{\ot^2}\ot \Pi^{\hs L}).
$$
\end{proposition}

\begin{proof} We have
$$

$$
where the first equality holds by Remark~\ref{rem1}; the second one, since $u_2*f = f$; the third one, since $\mu$ is associative and $c$ is natural; the fourth one, by Proposition~\ref{Delta compuesto con Pi^R multiplicado por algo}; and the last one, by Remark~\ref{epsilon mult * f=f}. Consequently $u_3\xst F_{\!\epsilon}\!=\! f\xcirc (H\ot \mu)\xcirc (H^{\ot^2}\ot \Pi^{\hs L})$. A similar argument proves that
$$
F_{\!\epsilon}\xst u_3 = f\xcirc (H\ot \mu)\xcirc (H^{\ot^2}\ot \ov{\Pi}^{\hs L})
$$
By Proposition~\ref{vale} this finishes the proof.
\end{proof}

\begin{definition}\label{regular map} A map $g\colon H^{\ot^n}\to A$ is {\em regular} if $g \xst  u_n = g$ and there exists $g'\in \Hom(H^{\ot^n},A)$ such that $g\xst  g' = g' \xst  g = u_n$.
\end{definition}

We let $\Reg(H^{\ot^n},A)$ denote the subset of $\Hom(H^{\ot^n},A)$ consisting of the regular maps. It is clear that $\Reg(H^{\ot^n},A)$ is closed under convolution product. Moreover, by Remark~\ref{u_n es idempotente} and Proposition~\ref{g*u_n} we can assume that $g' = g' \xst  u_n$ (and we do it). So, $\Reg(H^{\ot^n},A)$ is a group with identity element $u_n$ and $g'$ is the inverse of $g$. In the sequel we will write $g^{-1}$ instead of $g'$.

\begin{definition}\label{def cociclo inversible}
We say that the cocycle $f$ is {\em invertible} if it is regular.
\end{definition}

\begin{proposition}\label{alg prop inv} Let $F_1$, $F_2$ and $F_{\!\rho}$ be as in Proposition~\ref{alg prop}. If $f$ is invertible, then $F_1$, $F_2$ and $F_{\!\rho}$ are regular maps. Moreover
$$
F^{-1}_1 = f^{-1}\xcirc (\mu\ot H),\quad F^{-1}_2 = f^{-1}\xcirc (H\ot \mu)\quad\text{and}\quad F^{-1}_{\!\rho} = \rho \xcirc  (H\ot f^{-1}).
$$
\end{proposition}

\begin{proof} By Proposition~\ref{alg prop}(1) we know that $F_1 \xst u_3=u_3$. Moreover, since $\Delta$ is multiplicative and $f^{-1}\xst f = u_2$, we have
$$
(f^{-1}\xcirc (\mu\ot H)) \xst  F_1 = (f^{-1} \xst f)\xcirc (\mu\ot H)=u_2\xcirc (\mu\ot H) = u_3.
$$
Similarly, $F_1 \xst (f^{-1}\xcirc (\mu\ot H)) = u_3$. So, $F_1\in \Reg(H^{\ot^3},A)$ and $F_1^{-1} = f^{-1}\xcirc (\mu\ot H)$. A similar argument proves that $F_2$ is regular and $F_2^{-1} = f^{-1}\xcirc (H\ot \mu)$. Next we prove the assertion about~$F_{\!\rho}$. The fact that  $F_{\!\rho}\xst u_3= F_{\!\rho}$ it follows from Proposition~\ref{alg prop}(3). By the equalities in items~(2) and~(3) of Proposition~\ref{wbialgebras'} and the fact that $f^{-1}\xst f = u_2$, we have
$$
(\rho\xcirc (H\ot f^{-1})) \xst  F_{\!\rho}  = \rho\xcirc (H\ot f^{-1}\xst f)= \rho\xcirc (H\ot u_2)=u_3.
$$
Similarly $F_{\!\rho} \xst (\rho\xcirc (H\ot f^{-1})) = u_3$. So $F_{\!\rho}\in \Reg(H^{\ot^3},A)$ and $F^{-1}_{\!\rho} = \rho\xcirc (H\ot f^{-1})$.
\end{proof}

\begin{proposition}\label{alg prop inv'} Let $F_{\!\epsilon}$ be as in Proposition~\ref{F_epsilon*u_3 = u_3*F_epsilon} and let $F'_{\!\epsilon}\coloneqq f^{-1}\ot \epsilon$. If $f$ is invertible, then the equalities $F_{\!\epsilon}\xst F'_{\!\epsilon} = F'_{\!\epsilon} \xst F_{\!\epsilon} = u_2\ot\epsilon$ hold.
\end{proposition}

\begin{proof}
Since $f\xst f^{-1} = u_2$, we have $F_{\!\epsilon}\xst F'_{\!\epsilon} = (f\xst f^{-1})\ot \epsilon = u_2\ot \epsilon$. A similar argument,  using that $f^{-1}\xst f = u_2$, shows that $F'_{\!\epsilon}\xst F_{\!\epsilon} = u_2\ot\epsilon$.
\end{proof}

\begin{remark}\label{bar F epsilon} Let $F_{\!\epsilon}$ be as in Proposition~\ref{F_epsilon*u_3 = u_3*F_epsilon} and let $F'_{\!\epsilon}$ be as in Proposition~\ref{alg prop inv'}. Assume that $f$ is invertible. By Remark~\ref{u_n es idempotente} and Propositions~\ref{alg prop1}, \ref{F_epsilon*u_3 = u_3*F_epsilon} and~\ref{alg prop inv'},
$$
F'_{\!\epsilon}\xst u_3 \xst \hat{F}_{\!\epsilon} = F'_{\!\epsilon}\xst u_3 \xst u_3\xst F_{\!\epsilon} = F'_{\!\epsilon}\xst u_3\xst F_{\!\epsilon} = F'_{\!\epsilon}\xst F_{\!\epsilon}\xst u_3 = (u_2\ot\epsilon) \xst u_3 = u_3
$$
and
$$
\hat{F}_{\!\epsilon} \xst F'_{\!\epsilon}\xst u_3 = u_3\xst F_{\!\epsilon}\xst F'_{\!\epsilon}\xst u_3 = u_3\xst (u_2\ot\epsilon) \xst u_3 = u_3 \xst u_3 = u_3,
$$
where $\hat{F}_{\!\epsilon}\coloneqq u_3\xst F_{\!\epsilon}$. Consequently, $\hat{F}_{\!\epsilon} \in \Reg(H^{\ot^3},A)$ and $\hat{F}_{\!\epsilon}^{-1} = F'_{\!\epsilon}\xst u_3$.
\end{remark}

\begin{proposition}\label{equiv cond cocyclo} Let $F_1$, $F_2$ and $F_{\!\rho}$ be as in Proposition~\ref{alg prop} and let~$\hat{F}_{\!\epsilon}$ be as in Remark~\ref{bar F epsilon}. If $f$ is invertible, then $F_1,F_2,F_{\!\rho}, \hat{F}_{\!\epsilon} \in \Reg(H^{\ot^3},A)$ and $F_2\xst F_1^{-1} = F^{-1}_{\!\rho}\xst \hat{F}_{\!\epsilon}$.
\end{proposition}

\begin{proof}
The cocycle condition reads $F_{\!\rho}\xst F_2 = F_{\!\epsilon} \xst F_1$, where $F_{\!\epsilon}$ is as in Proposition~\ref{F_epsilon*u_3 = u_3*F_epsilon}. Since $F_{\!\epsilon} \xst F_1 = F_{\!\epsilon} \xst u_3 \xst F_1 = \hat{F}_{\!\epsilon} \xst F_1$, this implies that $F_{\!\rho}\xst F_2 = \hat{F}_{\!\epsilon} \xst F_1$. This finishes the proof because, by Proposition~\ref{alg prop inv} and the previous remark, $F_1,F_2,F_{\!\rho},\hat{F}_{\!\epsilon}\in \Reg(H^{\ot^3},A)$.
\end{proof}

\begin{proposition}\label{f^{-1} es normal}
Assume that $f$ is invertible. If $f$ is normal, then $f^{-1}$ is also.
\end{proposition}

\begin{proof} in fact,
$$
f^{-1}\xcirc (H\ot \eta) = ((f\xcirc (\mu\ot \eta))\xst f^{-1}) \xcirc (H\ot \eta) = (f\xst f^{-1})\xcirc (H\ot \eta) = u_2\xcirc (H\ot \eta) = u_1,
$$
where the first equality holds since $f^{-1}= u_2*f^{-1}$ and $u_2=u_1\xcirc \mu=f\xcirc (\mu\ot\eta)$; the second one, by Proposition~\ref{fundamental'}; and the third one, since $f\xst f^{-1}=u_2$. Similarly, $f^{-1}\xcirc (\eta\ot H)=u_1$.
\end{proof}

\begin{proposition}\label{fundamental''} Assume that $f$ is invertible. Then
\begin{align*}
&f^{-1}\xcirc (\mu\ot H)\xcirc (H\ot \Pi^{\hs R}\ot H) = f^{-1}\xcirc (H\ot \mu)\xcirc (H\ot \Pi^{\hs R}\ot H)
\shortintertext{and}
&f^{-1}\xcirc (\mu\ot H)\xcirc (H\ot \Pi^{\hs L}\ot H) = f^{-1}\xcirc (H\ot \mu)\xcirc (H\ot \Pi^{\hs L}\ot H).
\end{align*}
\end{proposition}

\begin{proof} By Remark~\ref{epsilon mult * f=f} and Propositions~\ref{equivalencia de que F incluido en A times V}, \ref{fundamental'} and~\ref{caso PiL}.
\end{proof}

\begin{definition}\label{inversa de gamma} When $f$ is invertible we define $\gamma^{-1}\colon H\to E$ by
$$
\gamma^{-1}\coloneqq \mu_E\xcirc (\jmath_{\nu} \ot \gamma)\xcirc Q,
$$
where $Q\coloneqq (f^{-1}\ot H)\xcirc(H\ot c)\xcirc (\Delta\ot H)\xcirc(S\ot H)\xcirc\Delta$.
\end{definition}

\begin{remark}\label{expresion alternativa de Q} Set $\bar{f}\coloneqq f^{-1}$. Note that
$$

$$
where $L\coloneqq f^{-1}\xcirc (S\ot H)\xcirc \Delta$.
\end{remark}

\begin{proposition}\label{gamma Pi^R * gamma^{-1}} Assume that $f$ is invertible. Then $(\gamma\xcirc \Pi^{\hs R})*\gamma^{-1} = \gamma^{-1}$.
\end{proposition}

\begin{proof} For legibility in the diagrams we set $\bar{\gamma}\coloneqq \gamma^{-1}$. We have
$$

$$
where the first equality holds by the definition of $\Pi^{\hs R}$; the second one, by  Propositions~\ref{Delta compuesto con unidad incluido en H^{hs R}ot H^{hs L}}(1) and~\ref{Delta compuesto con Pi^R multiplicado por algo}; the third one, by Definition~\ref{inversa de gamma}; the fourth one, by Propositions~\ref{Delta compuesto con unidad incluido en H^{hs R}ot H^{hs L}}(1) and~\ref{Delta compuesto con Pi^R multiplicado por algo}; the fifth one, since $S$ is antimultiplicative; the sixth one, by Proposition~\ref{auxiliar 3}; the seventh one, since $\gamma\xcirc \eta =\eta_E$; and the last one, by Definition~\ref{inversa de gamma} and Remark~\ref{expresion alternativa de Q}.
\end{proof}

\begin{remark}\label{remm1} Assume that $f$ is invertible and set $\bar{f}\coloneqq f^{-1}$ and $\bar{F}_{\rho}\coloneqq F_{\rho}^{-1}$. We have
$$

$$
where the first equality holds by Proposition~\ref{F_epsilon*u_3 = u_3*F_epsilon}; the second one, by the definition of $\chi_{\rho}$; the third one, by Proposition~\ref{delta Pi^R, etc}; the fourth one, by Proposition~\ref{fundamental''}; the fifth one, by Definition~\ref{weak bialgebra}(1) and the fact that $\mu$ is unitary; and the last one, by the definition of $Q$.
\end{remark}

\begin{lemma}\label{inv implica cleft p}
We have $(F_2 \xst \bar{F}_1)\circ (H\ot S\ot H)\circ (\mu\ot \Delta)\circ (\Pi^{R}\ot \Delta)  = \rho\circ (\Pi^{R}\ot u_1)$ .
\end{lemma}

\begin{proof} We compute
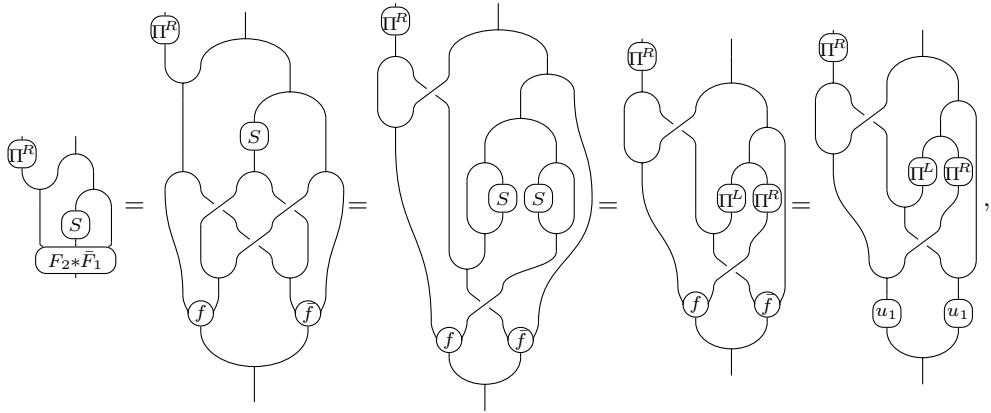
\begin{equation*}

\end{equation*}
where the first equality holds by the very definition of $F_2 \xst \bar{F}_1$; the second one, by Definition~\ref{weak bialgebra}, the discussion below Definition~\ref{weak Hopf Algebra} and the fact that $\Delta$ is coassociative and $\mu$ is associative; the third one, by Definition~\ref{weak Hopf Algebra}, and the last one, by by the fact that $f$ is normal, and Propositions~\ref{Delta compuesto con unidad incluido en H^{hs R}ot H^{hs L}}(1),~\ref{Delta compuesto con Pi^R}(4), \ref{fundamental'}, \ref{caso PiL}, \ref{f^{-1} es normal} and \ref{epsilon mult * f=f}. So, we are reduce to prove that the last diagram represent the function at the right hand of the equality in the statement. But this is true, since
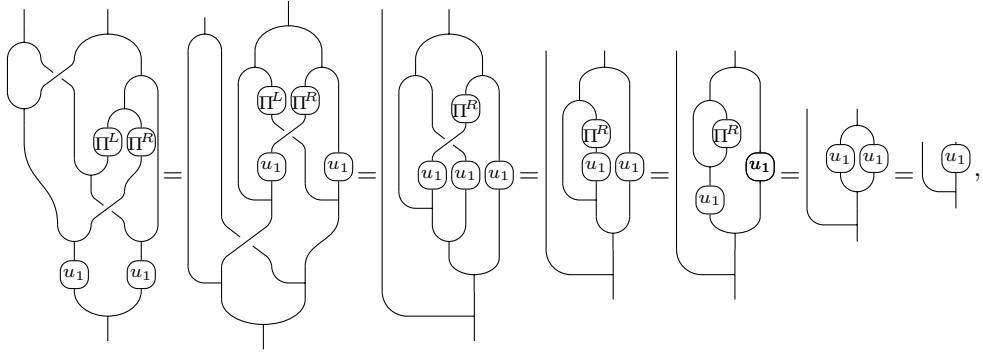
\begin{equation*}

\end{equation*}
where the first equality holds by the equality in Proposition~\ref{wbialgebras'}(3) and the coassociativity of $\Delta$; the second one, by the equalities in items~(2) and~(4) of Proposition~\ref{wbialgebras'} and the associativity of $\mu$; the third one, by the equality in  Proposition~\ref{wbialgebras'}(2) and the fact that $\mu$ is unitary; the fourth one, by the equality in  Proposition~\ref{wbialgebras'}(3); the fifth one, by Remark~\ref{ide ast Pi^R = ide and}; and the sixth one, since $u_1$ is idempotent.
\end{proof}

\begin{proposition}\label{inv implica cleft} Let $\gamma$ be as above of Theorem~\ref{prodcruz1}. If $f$ is invertible, then $\gamma^{-1}\xst \gamma = \gamma\xcirc \Pi^{\hs R}$ and $\gamma\xst \gamma^{-1} = \gamma\xcirc \Pi^{\hs L}$.
\end{proposition}

\begin{proof} For the legibility in the diagrams we set $\bar{f}\coloneqq\check{} f^{-1}$, $\bar{\gamma}\coloneqq \gamma$, $\bar{F}_1\coloneqq F_1^{-1}$ and $\bar{F}_{\!\rho}\coloneqq F_{\!\rho}^{-1}$. First we prove that $\bar{\gamma}\xst \gamma = \gamma\xcirc \Pi^{\hs R}$. Let $Q$ be as in Definition~\ref{inversa de gamma}. By Remark~\ref{aux}(3), and the fact that $\bar{f}\xst f=u_2$, $\Delta$ is multiplicative and $S\xst \ide=\Pi^{\hs R}$, we have
$$

$$
Using this, the associativity of $\mu_E$, Remark~\ref{gama iota y gama gama} and Theorem~\ref{weak crossed prod}(11), we obtain
$$
%
$$
as desired. We next prove that $\gamma\xst \bar{\gamma} = \gamma\xcirc \Pi^{\hs L}$. By Remark~\ref{gama iota y gama gama}, Theorem~\ref{weak crossed prod}(11) and the associativity of $\mu_E$,
$$
%
$$
Moreover, we have
$$
%
$$
where the first equality holds by the definition of $\mathcal{F}_{\hs f}$; the second one, by Remark~\ref{aux}(1) and the fact that, by Remark~\ref{expresion alternativa de Q},
$$
%
$$
the third one, since $\Delta$ is coassociative and $c$ is natural; the fourth one, by Remark~\ref{remm1}, the equality in Definiton~\ref{weak Hopf Algebra}(1), the definition of $\Pi^{\hs L}$ and the fact that $\Delta$ is coassociative; the fifth one, by Proposition~\ref{Delta compuesto con Pi^R multiplicado por algo} and the fact that $\Delta$ is counitary; and the sixth one, by Proposition~\ref{equiv cond cocyclo}. Consequently,
$$
%
$$
where the third equality holds by  Lemma~\ref{inv implica cleft p}, the fourth one, by the fact that $c$ is natural and the equalities in items~(3) and~(6) of Proposition~\ref{wbialgebras'}; the fifth one, by Proposition~\ref{Delta compuesto con Pi^R multiplicado por algo}; the sixth one, since $\mu$  is unitary; the seventh one, by the definitions of $u_1$ and $\chi_{\rho}$; the eighth one, by  Remark~\ref{gama iota y gama gama}; and the last one, by Theorem~\ref{weak crossed prod}(11).
\end{proof}

\begin{proposition}\label{es cleft} The equality $\delta_E\xcirc \gamma\xcirc \Pi^{\hs L}=(E\ot \Pi^L)\xcirc \delta_E\xcirc \gamma\xcirc \Pi^{\hs L}$ holds.
\end{proposition}

\begin{proof} We have
$$
\delta_E\xcirc \gamma\xcirc \Pi^{\hs L} = (\gamma\ot H)\xcirc \Delta \xcirc \Pi^{\hs L} = (\gamma\ot \Pi^{\hs L})\xcirc \Delta \xcirc \Pi^{\hs L} = (H\ot\Pi^{\hs L})\xcirc \delta_E\xcirc \gamma\xcirc \Pi^{\hs L},
$$
where the first and last equalities hold since $\gamma$ is right $H$-colinear; and the second one, by Pro\-position~\ref{mu delta Pi^R, etc}
%
%
%
\end{proof}

\begin{theorem}\label{crossed prod cleft} Let $A$ be a weak $H$-module algebra with weak action $\rho$ and let $f\colon H\ot H\to A$ be a map. Assume that the hypotheses of Theorem~\ref{weak crossed prod} are fulfilled and let $E$ be the unitary crossed product of $A$ by $H$ associated with $\rho$ and $f$. If $f$ is convolution invertible, then $(E,\jmath_{\nu})$ is $H$-cleft.
\end{theorem}

\begin{proof} By Propositions~\ref{E es H-comodulo algebra debil}, \ref{el egalizador es A}, \ref{gamma Pi^R * gamma^{-1}}, \ref{inv implica cleft} and~\ref{es cleft}.
\end{proof}

\begin{remark}\label{prod cruzados equiv equivale a extensiones equiv}
Let $\cramped{E'=A\times_{\rho'}^{f'} H}$ be another unitary crossed product satisfying the hypotheses of Theorem~\ref{crossed prod cleft}. Let $\nu'$ be the preunit of $\cramped{A\ot_{\rho'}^{f'} H}$. By Proposition~\ref{el egalizador es A} and the fact that an unitary algebra morphism $\Phi\colon E\to E'$ is left $A$-linear if and only if $\Phi\xcirc \jmath_{\nu}= \jmath_{\nu'}$, the crossed products $E$ and $E'$ are equivalent if and only if the extensions $(E,\jmath_{\nu})$ and $(E',\jmath_{\nu'})$ are equivalent.
\end{remark}

\begin{proposition}\label{propiedad 4} The equality
$$
\mu_E\xcirc (\jmath_{\nu}\ot \gamma^{-1})\xcirc c = \mu_E\xcirc (\gamma^{-1}\ot \jmath_{\nu})\xcirc (H\ot \rho)\xcirc (\Delta\ot A)
$$
holds.
\end{proposition}

\begin{proof} By the legibility in the diagrams we set $\bar{\gamma}\coloneqq \gamma^{-1}$. We have
$$

$$
where the first equality follows by the equality in Proposition~\ref{wbialgebras'}(2); the second one, from Propositions~\ref{auxiliar4} and~\ref{inv implica cleft}; the third one, by Remark~\ref{aux}(1) and the associativity of $\mu_E$; and the last one, by the first equality in Remark~\ref{gama iota y gama gama}. Using this and Propositions~\ref{conmutatividad de gamma e i}, \ref{gamma Pi^R * gamma^{-1}} and~\ref{inv implica cleft}, we obtain
$$
%
$$
as desired.
\end{proof}

\section{Cleft extensions are weak crossed products with invertible cocycle}\label{Section: Cleft extensions are weak crossed products with invertible cocycle}

Let $H$ be a weak bialgebra. In this Section we prove that each $H$-cleft extension is isomorphic to a crossed product with invertible cocycle, of $H$ by a weak $H$-module algebra. In a final remark we prove that when $H$ is a weak Hopf algebra, the category of unitary crossed products of $A$ by $H$ with invertible cocycle, such that $A$ a weak $H$-module algebra, and the category of $H$-cleft extensions of $A$ are equivalent.

\smallskip

Let $(E,\jmath)$ be a cleft extension of $A$ by $H$ and let $\gamma\colon H\to  E$ be a convolution invertible total integral. Let $\Upsilon\colon H\ot E\longrightarrow E\ot H$ be the map defined by
$$
\Upsilon\coloneqq (E\ot \mu)\xcirc (c\ot H)\xcirc (H\ot \delta_E).
$$
It is well know that $(E,H,\Upsilon)$ is a weak entwining data and $(E,\mu_E,\delta_E)$ is a weak entwined module (for the definitions of weak entwining data and weak entwined module see \cite{AFGR}*{Definitions~3.1 and~3.2}). Let $q_{\gamma^{-1}}\colon E\to E$ and $w_{\gamma}^E\colon A\ot H\to E$ be the maps defined by
$$
q_{\gamma^{-1}}^E\coloneqq \mu_E\xcirc (E\ot \gamma^{-1})\xcirc \delta_E\quad\text{and}\quad w_{\gamma}^E\coloneqq \mu_E \xcirc (\jmath\ot \gamma),
$$
respectively. By \cite{AFGR}*{3.8} there exists a morphism $p_{\gamma^{-1}}\colon E\to A$ such that $q_{\gamma^{-1}} = \jmath \circ p_{\gamma^{-1}}$. Let $\tilde{w}_{\gamma^{-1}}^E \colon E\to A\ot H$ be the morphism defined by $\tilde{w}_{\gamma^{-1}}^E\coloneqq (p_{\gamma^{-1}}^E\ot H)\xcirc \delta_E$. By \cite{AFGR}*{3.10} we know that $w_{\gamma}^E\xcirc \tilde{w}_{\gamma^{-1}}^E = \ide_E$, so that the map $\Omega_E\coloneqq \tilde{w}_{\gamma^{-1}}^E\xcirc w_{\gamma}^E$ is an idempotent. 

\begin{remark}\label{igualdad 0} Set $\bar{\gamma}\coloneqq \gamma^{-1}$. Since $\gamma$ is right $H$-colinear, we have
$$

$$
where the first and fourth equality holds since $(E,H,\Upsilon)$ is an entwined module; the second one, by the equality in Proposition~\ref{wbialgebras}(3) and the fact that, since by Proposition~\ref {algunas composiciones} and the fact $\jmath$ is an extension,
$$
\delta_E \xcirc \jmath = (E\ot \Pi^L)\xcirc \delta_E \xcirc \jmath = (E\ot\ov{\Pi}^R)\xcirc \delta_E \xcirc \jmath;
$$
the third one, since $\mu_E$ is associative; and the last one,
since $\mu_E$ is unitary.
\end{remark}

\begin{remark}\label{w_E es lineal a izquierda} Set $w\coloneqq w_{\gamma}^E$. We have
$$
%
$$
where the first, fourth, fifth and last equalities hold by the definition of $w$; the second one, since $\jmath$ is multiplicative, the third one, since $\mu_E$ is associative; the sixth one, by Remark~\ref{igualdad 1}; and the seventh one, since $\gamma$ is right $H$-colinear. Consequently $w$ is left $A$-linear and right $H$-colinear.
\end{remark}

\begin{proposition}\label{q y multiplicacion nivel 1} The equality $p_{\gamma^{-1}}^E\xcirc \mu_E \xcirc (\jmath \ot E)= \mu_A\xcirc (A \ot p_{\gamma^{-1}}^E)$ holds.
\end{proposition}

\begin{proof} Set $\bar{\gamma}\coloneqq \gamma^{-1}$, $q\coloneqq q_{\gamma^{-1}}^E$ and $p\coloneqq p_{\gamma^{-1}}^E$. We have
$$

$$
where the first and fourth equality are true by the very definition of $q$; the second one, by Re\-mark~\ref{igualdad 1}; the third one, by the associativity of $\mu_E$; and the last one, because $q=\jmath\xcirc p$ and $\jmath$ is an algebra morphism. Since $q=\jmath\xcirc p$ and $\jmath$ is a monomorphism, this finishes the proof.
\end{proof}

\begin{remark}\label{w'E es E_C lineal a izq} Set $p\coloneqq p_{\gamma^{-1}}^E$ and $\tilde{w}\coloneqq \tilde{w}_{\gamma^{-1}}^M$. We have
$$
%
$$
where the first, fourth, fifth and the last equalities hold by the very definition of $\tilde{w}$; the second one, by Remark~\ref{igualdad 1}; the third one, by Proposition~\ref{q y multiplicacion nivel 1}; and the sixth one, since $\delta_E$ is coassociative. Consequently, $\tilde{w}$ is left $A$ liner and right $H$-colinear.
\end{remark}


\begin{proposition}\label{multiplicacon en Ec ot C es Ec lineal a izq}
The map $\tilde{\mu}\colon A\ot H \ot A \ot H\longrightarrow A\ot H$, defined by $\tilde{\mu}\coloneqq \tilde{w}_{\gamma^{-1}}^E\xcirc \mu_E\xcirc (w_{\gamma}^E\ot w_{\gamma}^E)$, is an associative product which is normalized with respect to $\Omega_E$. Moreover $\tilde{\mu}$ is left $A$-linear, $\tilde{\nu}\coloneqq \tilde{w}_{\gamma{-1}}^E\xcirc \eta_E$ is a preunit of $\tilde{\mu}$ and $\Omega_E = \nabla_{\!\tilde{\nu}}$.
\end{proposition}

\begin{proof} By Remark~\ref{producto asociado con preunidad} we know that $\tilde{\mu}$ is an associative product that is normalized respect to~$\Omega_E$, that $\tilde{\nu}$ is a preunit of $\tilde{\mu}$ and that $\Omega_E = \nabla_{\!\tilde{\nu}}$. Set $w\coloneqq w_{\gamma}^E$ and $\tilde{w}\coloneqq \tilde{w}_{\gamma^{-1}}^M$. Since, by the associativity of $\mu_E$ and Remarks~\ref{w_E es lineal a izquierda} and~\ref{w'E es E_C lineal a izq},
$$

$$
the morphism $\tilde{\mu}$ is left $A$-linear. We next prove that $\tilde{\mu}$ is right $H$-colinear. By Remarks~\ref{w_E es lineal a izquierda} and~\ref{w'E es E_C lineal a izq} we know that $w$ and $\tilde{w}$ are right $H$-colinear. Since $\mu_E$ is also right $H$-colinear, we have
$$
%
$$
as desired.
\end{proof}

\begin{theorem}\label{teo 1} Let $\tilde{\mu}$ and $\tilde{\nu}$ be as in Proposition~\ref{multiplicacon en Ec ot C es Ec lineal a izq}. The morphisms
$$
\rho\colon H\ot A\longrightarrow A\qquad\text{and}\qquad f\colon H\ot H\longrightarrow A,
$$
defined by
$$
\rho\coloneqq (A\ot\epsilon)\xcirc \tilde{\mu}\xcirc \bigl(\eta_A\ot H\ot \jmath'_{\tilde{\nu}}\bigr)\quad\text{and}\quad f\coloneqq (A\ot\epsilon)\xcirc \tilde{\mu}\xcirc \bigl(\eta_A\ot H\ot \eta_A\ot H\bigr),
$$
where $\jmath'_{\tilde{\nu}}\coloneqq (\mu_A\ot H)\xcirc (A\ot \tilde{\nu})$, satisfy the following properties:

\begin{enumerate}[itemsep=0.7ex,  topsep=1.0ex, label=\emph{(\arabic*)}]

\item $\rho$ is a weak measure of $H$ on $A$,

\item $f$ is a cocycle that satisfies the twisted module condition,

\item $f = \mu_A\xcirc(A \ot\rho)\xcirc (f\ot \mu\ot A)\xcirc (\Delta_{H\ot H}\ot\eta_A)$,

\item $\rho\xcirc (H\ot \eta_A) = \mu_A\xcirc (\rho\ot f) \xcirc (H\ot c\ot H)\xcirc (\Delta\ot \tilde{\nu})$,

\item $\rho\xcirc (H\ot \eta_A) = \mu_A\xcirc (A\ot f)\xcirc (\tilde{\nu}\ot H)$,

\item $(\mu_A\ot H)\xcirc (A\ot\rho\ot H)\xcirc (A\ot H\ot c)\xcirc (A\ot \Delta\ot A)\xcirc(\tilde{\nu}\ot A) = (\mu_A\ot H)\xcirc (A\ot \tilde{\nu})$.

\end{enumerate}
Moreover
$$
\tilde{\mu}=\mu_{A\ot_{\rho}^f H},\qquad E \simeq A\times_{\rho}^f H \qquad \text{and} \qquad \gamma = w_{\gamma} \xcirc (\eta_A\ot H).
$$
\end{theorem}

\begin{proof} Since $\tilde{\nu} = \tilde{w}_{\gamma^{-1}}^E \xcirc \eta_E = (p_{\gamma^{-1}}^E\ot H)\xcirc \delta_E \xcirc \eta_E$, from the equality in Proposition~\ref{wbialgebras}(6), it follows that $\tilde{\nu} = (A\ot\Pi^{\hs L})\xcirc \tilde{\nu}$. Thus, by \cite{FGR}*{Theorems~3.11 and~4.2}, in order to prove the result it suffices to note that, by Proposition~\ref{multiplicacon en Ec ot C es Ec lineal a izq}, the map $\tilde{\mu}$ is left $A$-linear, right $H$-colinear and associative, the map $\tilde{\nu}$ is a preunit of $\tilde{\mu}$; and $\tilde{\mu}$ is normalized with respect to $\nabla_{\!\tilde{\nu}}$.
\end{proof}

\begin{remark}\label{calculos} Set $w\coloneqq w_{\gamma}^E$, $\tilde{w}\coloneqq \tilde{w}_{\gamma^{-1}}^E$ and $p\coloneqq p_{\gamma^{-1}}^E$. By Remark~\ref{w'E es E_C lineal a izq},
$$

$$
Consequetly, by the definitions of $\tilde{\mu}$, $w$, $\tilde{w}$, $\rho$ and $f$, we have
$$
%
$$
Moreover $\jmath_{\tilde{\nu}} = w\xcirc \tilde{w}\xcirc \jmath = \jmath$.
\end{remark}

\begin{remark}\label{prev} By Remark~\ref{igualdad 0} and the fact that $\jmath \xcirc p_{\gamma^{-1}}^E = q_{\gamma^{-1}}^E$, $\gamma\xcirc \eta = \eta _E$ and $\gamma * \gamma^{-1} = \gamma\xcirc \Pi^L$, we have
$$
\jmath \xcirc p_{\gamma^{-1}}^E\xcirc \eta_E = q_{\gamma^{-1}}^E\xcirc \gamma\xcirc \eta = (\gamma * \gamma^{-1})\xcirc \eta = \gamma\xcirc \Pi^L\xcirc\eta = \gamma\xcirc \eta = \eta_E = \jmath\xcirc \eta_A.
$$
Since $\jmath$ is a monomorphism, this implies that $p_{\gamma^{-1}}^E\xcirc \eta_E = \eta_A$. Consequently
$$
p_{\gamma^{-1}}^E\xcirc \jmath = p_{\gamma^{-1}}^E\xcirc \mu_E \xcirc (\jmath \ot \eta_E) = \mu_A\xcirc (A \ot p_{\gamma^{-1}}^E) \xcirc (A \ot \eta_E) = \ide_A,
$$
where the second equality holds by Proposition~\ref{q y multiplicacion nivel 1}.
\end{remark}

\begin{remark}By Theorem~\ref{teo 1} the map $\gamma$ satisfies Lemma~\ref{caso particular} and Proposition~\ref{multiplicacion y gamma}. We will use these facts in the proof of Proposition~\ref{es un H modulo algebra debil} and Lemma~\ref{algo}.
\end{remark}

\begin{proposition}\label{es un H modulo algebra debil} The algebra $A$ is a left weak $H$-module algebra.
\end{proposition}

\begin{proof} Since $\rho$ is a weak measure we know that the equality in Proposition~\ref{wbialgebras'}(2) is satisfied. Thus, in order to finish the proof it suffices to check that the equalities in items~(1), (3) and~(7)~of Proposition~\ref{wbialgebras'} also are satisfied.
For the legibility in the diagrams we set $p\coloneqq p_{\gamma^{-1}}^E$ and $q\coloneqq q_{\gamma^{-1}}^E$. For the equality in Proposition~\ref{wbialgebras'}(1) we have
$$
\rho\xcirc (\eta\ot A) = p\xcirc \mu_E \xcirc (\gamma\ot \jmath) \xcirc (\eta\ot A) = p\xcirc \jmath = \ide_A,
$$
where the first equality holds by Remark~\ref{calculos}; the second one, since $\gamma\xcirc \eta = \eta_E$; and the last one, by Remark~\ref{prev}. We next prove that the equality in Proposition~\ref{wbialgebras'}(3) is true. We have
$$

$$
where the first and last equality hold since $\jmath$ and $\mu_E$ are unitary and $\jmath\xcirc p = q$; the second, fourth and sixth one, by Remark~\ref{igualdad 0} and the fact that $\gamma * \gamma^{-1} = \gamma\xcirc \Pi^L$; the third one, by Proposition~\ref{multiplicacion y gamma}(1); and the fifth one, by Proposition~\ref{mu delta Pi^R, etc}. This, combined with Remark~\ref{calculos} and the fact that $\jmath\xcirc p = q$ and $\jmath$ is a monomorphism, proves Proposition~\ref{wbialgebras'}(3). Finally, the equality in Proposition~\ref{wbialgebras'}(6) is satisfied, since
$$

$$
where the first and last equalities hold since $\jmath$ and $\mu_E$ are unitary; the second one, by Remark~\ref{igualdad 0} and the fact that $\gamma * \gamma^{-1} = \gamma\xcirc \Pi^L$; and the third one, since $p\xcirc q = q$.
\end{proof}

Next we are going to prove that $f$ is regular. By Theorem~\ref{teo 1}(3) we know that $f\xst u_2 = f$. Moreover, by Remark~\ref{a derecha implica aizquierda} and the equality in Proposition~\ref{wbialgebras'}(3), we also have $u_2\xst f = f$. Let $\sigma_E\colon H\ot H\to E$ be the morphism defined by
$$
\sigma_E\coloneqq \mu_E\xcirc (\mu_E\ot \gamma^{-1}) \xcirc (\gamma\ot \Upsilon)\xcirc (\Delta\ot \gamma).
$$
Clearly $\sigma_E = \bigl(\mu_E\xcirc (\gamma\ot \gamma)\bigr) \xst (\gamma^{-1}\xcirc \mu)$. By Remark~\ref{calculos} and \cite{AFGR2}*{Proposition~1.17} we know that
$$
u_2 = p_{\gamma^{-1}}\xcirc \gamma\xcirc \mu\quad\text{and}\quad \sigma_E = \jmath \xcirc f.
$$
Note that
\begin{align*}
&(q_{\gamma^{-1}}\xcirc \gamma\xcirc \mu)\xst \sigma_E = (\jmath\xcirc u_2)\xst (\jmath\xcirc f) = \jmath\xcirc (u_2\xst f) = \jmath\xcirc f = \sigma_E
\shortintertext{and}
&\sigma_E\xst (q_{\gamma^{-1}}\xcirc \gamma\xcirc \mu) =  (\jmath\xcirc f)\xst (\jmath\xcirc u_2) = \jmath \xcirc (f\xst u_2) = \jmath\xcirc f = \sigma_E.
\end{align*}

\begin{remark}\label{para usar} Since $\sigma_E$ and $q_{\gamma^{-1}}\xcirc \gamma$ factorize through $\jmath$,
$$
\delta_E\xcirc \sigma_E = (E\ot \Pi^{\hs L})\xcirc \delta_E\xcirc \sigma_E\quad\text{and}\quad \delta_E\xcirc q_{\gamma^{-1}}\xcirc \gamma= (E\ot \Pi^{\hs L})\xcirc \delta_E\xcirc q_{\gamma^{-1}}\xcirc \gamma.
$$
\end{remark}

\begin{lemma}\label{algo} Let $\sigma_E^{-1}\colon H\ot H\to E$ be the map defined by $\sigma_E^{-1}\coloneqq (\gamma \xcirc \mu)\xst \bigl(\mu_E\xcirc c\xcirc (\gamma^{-1}\ot \gamma^{-1})\bigr)$. The following equalities hold:
\begin{equation}\label{A1}
\sigma_E\xst \sigma_E^{-1} = q_{\gamma^{-1}}\xcirc \gamma\xcirc \mu\qquad\text{and}\qquad \sigma_E^{-1}\xst \sigma_E = q_{\gamma^{-1}}\xcirc \gamma\xcirc \mu
\end{equation}
\end{lemma}

\begin{proof} First we show that the first equality in~\eqref{A1} is satisfied. Since $\mu$ is comultiplicative,
$$
(\gamma^{-1}\xcirc \mu)\xst (\gamma \xcirc \mu) = (\gamma^{-1}\xst \gamma) \xcirc \mu = \gamma\xcirc\Pi^{\hs R}\xcirc \mu.
$$
We claim that
$$
\bigl(\mu_E\xcirc (\gamma\ot \gamma)\bigr) \xst (\gamma\xcirc\Pi^{\hs R}\xcirc \mu) = \mu_E\xcirc (\gamma\ot \gamma).
$$
In fact, this is true since
$$

$$
where the first equality holds since $\mu_E$ is associative; the second one, by Proposition~\ref{multiplicacion y gamma}(2); the third one, by Proposition~\ref{Delta compuesto con Pi^R}(2); the fourth one, by Remark~\ref{ide ast Pi^R = ide and} and the fact that $\Delta$ is coassociative and $c$ is natural; the fifth one, by Proposition~\ref{mu Pi^R, etc}(2); the sixth one, by Proposition~\ref{Delta eta con S}(1); and the last one, by Lemma~\ref{caso particular}. For the sake of legibility in the following diagramas we set $\bar{\gamma}\coloneqq \gamma^{-1}$ and $q\coloneqq q_{\gamma^{-1}}^E$. In order to end the proof of the first equality in~\eqref{A1} it suffices to note that
$$

$$
where the first equality holds since $c$ is natural and $\mu_E$ is associative; the second and sixth one, since $\gamma\xst \bar{\gamma} = \gamma\xcirc \Pi^L$; the third one, by Proposition~\ref{multiplicacion y gamma}(1); the fourth and seventh one, by Proposition~\ref{mu Pi^R, etc}(1); the fifth one, since $\Delta$ is coassociative and $c$ is natural; the eighth one, by Proposition~\ref{mu delta Pi^R, etc}; and the last one, by Remark~\ref{igualdad 0} and the fact that $\gamma\xcirc \Pi^L = \gamma\xst\bar{\gamma}$. We next prove the second equality in~\eqref{A1}. To begin with we have
$$

$$
where the first equality holds since $c$ is natural and $\mu_E$ is associative; the second and sixth one, since $\bar{\gamma}\xst \gamma = \gamma\xcirc \Pi^{\hs R}$; the third one, by Proposition~\ref{multiplicacion y gamma}(2); the fourth and seventh one, by Proposition~\ref{mu Pi^R, etc}(2); the fifth one, since $\Delta$ is coassociative and $c$ is natural; and the last one, by Proposition~\ref{mu delta Pi^R, etc}. Thus the proof of the second equality in~\eqref{A1} follows, because
$$
%
$$
where the first and third equalities hold by Definition~\ref{weak bialgebra}(1); the second one, by Corollary~\ref{gamma xst (gamma xcirc Pi^R) = gamma}; and the last one, by Remark~\ref{igualdad 0}.
\end{proof}

\begin{remark}\label{esta donde debe} We have
$$
\sigma_E^{-1}\xst (q_{\gamma^{-1}}\xcirc \gamma\xcirc \mu) = \sigma_E^{-1}\xst \sigma^E \xst \sigma_E^{-1} = (q_{\gamma^{-1}}\xcirc \gamma\xcirc \mu) \xst \sigma_E^{-1} = \sigma_E^{-1},
$$
where the last equality follows from the definition of $\sigma_E^{-1}$ and the fact that, by~Remark~\ref{igualdad 0}, equality $\gamma \xst \gamma^{-1}=\gamma\circ \Pi^L$ and Corollary~\ref{gamma xst (gamma xcirc Pi^R) = gamma},
$$
(q_{\gamma^{-1}}\xcirc \gamma\xcirc \mu) \xst (\gamma \xcirc \mu) = \bigl((q_{\gamma^{-1}}\xcirc \gamma) \xst \gamma \bigr)\xcirc \mu = \bigl((\gamma\xcirc \Pi^{\hs L}) \xst \gamma \bigr)\xcirc \mu =  \gamma \xcirc \mu.
$$
\end{remark}

\begin{remark}\label{use} Set $\bar{\sigma}_E\coloneqq \sigma_E^{-1}$. We have
$$

$$
where the first and fifth equality hold by Proposition~\ref{subalgebras}(1) and Remark~\ref{para usar}; the second one, since $\delta_E$ is multiplicative; the third one, by Lemma~\ref{algo}; the fourth one, by Remark~\ref{para usar}; and the last one, by Remark~\ref{esta donde debe} and the fact that $\delta_E$ is multiplicative.
\end{remark}

\begin{theorem}\label{se factoriza} The cocycle $f$ is invertible.
\end{theorem}

\begin{proof} In order to abbreviate expressions we set $U\! =\! \delta_E \xcirc \sigma_E$, $\bar{U}\! =\! \delta_E \xcirc \sigma_E^{-1}$ and $N\! =\! \delta_E\xcirc q_{\gamma^{-1}}^E\xcirc \gamma\xcirc \mu$. We have
$$
\bar{U} = N\xst \bar{U} = ((E\ot \Pi^{\hs L})\xcirc \bar{U})\xst U\xst \bar{U} = ((E\ot \Pi^{\hs L})\xcirc \bar{U})\xst N = (E\ot \Pi^{\hs L})\xcirc \bar{U},
$$
where the first equality holds by Remark~\ref{esta donde debe} and the fact that $\delta_E$ is multiplicative; the second one, by the first part of Remark~\ref{use}; the third one, by Lemma~\ref{algo} and the fact that $\delta_E$ is multiplicative; and the last one, by the second part of Remark~\ref{use}. Consequently $\sigma^{-1}_E$ factorize through $\jmath$. Let $f^{-1}\colon H\ot H\to A$ be such that $\sigma^{-1}_E = \jmath\xcirc f^{-1}$. Since $\jmath$ is a monomorphism and
\begin{align*}
& \jmath\xcirc (f\xst f^{-1}) = (\jmath\xcirc f)\xst (\jmath\xcirc f^{-1}) = \sigma_E \xst \sigma_E^{-1} = q_{\gamma^{-1}}\xcirc \gamma\xcirc \mu = \jmath \xcirc p_{\gamma^{-1}}\xcirc \gamma\xcirc \mu
\shortintertext{and}
& \jmath\xcirc (f^{-1}\xst f) = (\jmath\xcirc f^{-1})\xst (\jmath\xcirc f) = \sigma_E^{-1} \xst \sigma_E = q_{\gamma^{-1}}\xcirc \gamma \xcirc \mu = \jmath \xcirc p_{\gamma^{-1}}\xcirc \gamma\xcirc \mu,
\end{align*}
we obtain that $f\xst f^{-1} = p_{\gamma^{-1}}\xcirc \gamma\xcirc \mu$ and $f^{-1} \xst f = p_{\gamma^{-1}}\xcirc \gamma\xcirc \mu$, as desired.
\end{proof}

\begin{theorem}\label{cleft implica crossed product con cociclo inversible}
Let $(E,\jmath)$ be a $H$-cleft extension of $A$ by $H$ and let $\rho$, $f$ and $\tilde{\nu}$ be as in Theorem~\ref{teo 1}. Then $A$ is a weak $H$-module algebra via $\rho$, the hypotheses of Theorem~\ref{weak crossed prod} are fulfilled, the cocycle $f$ is invertible and  $E\simeq A\times_{\rho}^f H$.
\end{theorem}

\begin{proof}
By Theorem~\ref{teo 1}, Proposition~\ref{es un H modulo algebra debil} and Theorem~\ref{se factoriza}.
\end{proof}

\begin{remark} From Theorem~\ref{crossed prod cleft}, Remark~\ref{prod cruzados equiv equivale a extensiones equiv} and Theorem~\ref{cleft implica crossed product con cociclo inversible} it follows that if $H$ is a weak Hopf algebra and $A$ is an algebra, then the category of unitary crossed products of $A$ by $H$ with invertible cocycle and $A$ a weak $H$-module algebra, and the category of $H$-cleft extensions of $A$, are equivalent.
\end{remark}

\end{document}